\documentclass[10pt]{article}
\usepackage{anysize}
\marginsize{2.0cm}{2.0cm}{2.0 cm}{2.0 cm}
\usepackage{setspace}

\usepackage{latexsym,amsfonts,amsmath}
\usepackage{graphicx}
\usepackage{color}
\usepackage{epstopdf}
\usepackage[titletoc]{appendix}

\newtheorem{condition}{Condition}[section]{\bfseries}{\itshape}

{\bfseries}{\itshape}

\newtheorem{theorem}{Theorem}[section]{\bfseries}{\itshape}

\newtheorem{corollary}{Corollary}[section]{\bfseries}{\itshape}

\newtheorem{proposition}{Proposition}[section]{\bfseries}{\itshape}

\newtheorem{example}{Example}[section]{\bfseries}{\itshape}

\newtheorem{lemma}{Lemma}[section]{\bfseries}{\itshape}

\newtheorem{remark}{Remark}[section]{\bfseries}{\itshape}

\newtheorem{definition}{Definition}[section]{\bfseries}{\itshape}

{\bfseries}{\itshape}

\begin{document}

\title{Gradual-impulsive control for continuous-time Markov decision processes with total undiscounted costs and constraints: linear programming approach via a reduction method}

\date{}

\author{Alexey Piunovskiy\thanks{Department of Mathematical Sciences, University of
Liverpool, Liverpool, U.K.. E-mail: piunov@liv.ac.uk.}~ and Yi
Zhang \thanks{Corresponding author. School of Mathematics, University of Birmingham, Edgbaston,
Birmingham,
B15 2TT, U.K.. Email: y.zhang.29@bham.ac.uk; or yi.zhang@liv.ac.uk}}

\maketitle

\par\noindent{\bf Abstract:} We consider the constrained optimal control problem for the gradual-impulsive CTMDP model with the performance criteria being the expected total undiscounted costs (from the running cost and the cost from each time an impulse being applied). The  discounted model is covered as a special case. We justify fully a reduction method, and close an open issue in the previous literature. The reduction method induces an equivalent but simpler standard CTMDP model with gradual control only, based on which, we establish effectively, under rather natural conditions, a linear programming approach for solving the concerned constrained optimal control problem.
\bigskip

\bigskip

\par\noindent
{\bf AMS (2020) subject classification:} Primary 90C40,  Secondary
60J75, 49N25

\section{Introduction}

The present paper investigates continuous-time Markov decision processes (CTMDPs) in Borel state and action spaces, where the decision maker can control the process via its local characteristics (transition rate), and also can control directly the state of the process. Such a model is called the gradual-impulsive (control) model. For the gradual-impulsive CTMDP model, we are concerned with the following constrained optimal control problem: the expected total undiscounted cost is to be minimized, subject to other performance measures (objectives) in the same form not exceeding predetermined levels.

The gradual-impulsive CTMDP model is quite general. It has two important sub-models. One is the standard CTMDP model, in which the decision maker only controls the transition rate of the process. The other one is the impulsive control model, in which the decision maker can only control instantaneously the state of the process. Each of them has a vast literature: for standard CTMDP models, see the monographs \cite{Guo:2009,Kitaev:1995,Prieto-RumeauHernandez-Lerma:2012} and the more recent one \cite{PiunovskiyZhang:2020Book}, which is influenced by \cite{Feinberg:2004,Feinberg:2012OHL}; for the impulsive control model, see e.g., \cite{Costa:1989,deSaporta,Gatarek:1992,Presman:2006}. (The latter references actually dealt with a more general class of processes than what is of concern here, namely, piecewise deterministic processes, see \cite{Davis:1993}.) The optimal stopping problem is an important example of impulsive control models, where the decision maker can decide when to stop the process, applying the impulse once and for all, see e.g., \cite{BauerlePopp:2018,Costa:1988}. Compared to the aforementioned two sub-models, there is relatively less literature on gradual-impulsive CTMDP models, see e.g., \cite{DufourPiunovskiy:2016,DufourPiunovskiy:2016a,Plum:1991,vanderDuynSchouten:1983,Yushkevich:1983,Yushkevich:1988}.

Most of the previous literature on gradual-impulsive CTMDP models allows one to apply at most one impulse at a given time moment, and the effect of the impulse is often deterministic, as in the recent work \cite{Miller:2020}. In the gradual-impulsive CTMDP model considered in this paper an impulse can be applied at any time moment, and one can apply multiple impulses at a single time moment. Such gradual-impulsive CTMDP models were considered in \cite{Yushkevich:1983,Yushkevich:1988} and more recently in \cite{DufourPiunovskiy:2016,DufourPiunovskiy:2016a}. In \cite{Yushkevich:1988}, which is a refinement of \cite{Yushkevich:1983}, and in \cite{DufourPiunovskiy:2016,DufourPiunovskiy:2016a}, the authors handled multiple simultaneous impulses  by extending the time or the state suitably, and after that, developed a theory for the resulting new point process, which is more complicated than the uncontrolled version of the original process. The analysis in \cite{DufourPiunovskiy:2016,DufourPiunovskiy:2016a,Yushkevich:1983,Yushkevich:1988} as well as in \cite{Miller:2020} is direct, in the sense that no connection with the standard CTMDP model was explored therein. A different method (the so called time discretization method) was taken in \cite{Plum:1991,vanderDuynSchouten:1983}, where the gradual-impulsive CTMDP model and the associated optimal control problem were studied as the limit of a sequence of discrete-time problems (for the skeleton models). The skeleton models are with complicated transition probabilities.

The present paper differs from the previous literature in terms of the problem statement and the method of investigations. The concerned optimal control problem for the gradual-impulsive CTMDP considered in the majority of the previous literature is unconstrained (with a single objective), as is the case in  \cite{DufourPiunovskiy:2016,Miller:2020,Plum:1991,Yushkevich:1983,Yushkevich:1988} as well as in \cite{Costa:2000,DempsterYe:1995,GuoKurashimaPiunovskiyZhang:2021}. The main optimality result in these papers was the establishment of the Bellman (optimality) equation, which is used to characterize and show the existence of optimal strategies (often known as the dynamic programming method). One of the relatively few works dealing with constrained problems for gradual-impulsive CTMDP models is \cite{DufourPiunovskiy:2016a}, where the performance criteria are the expected total discounted costs, and a linear programming approach was established. The linear program formulation in \cite{DufourPiunovskiy:2016a} is a consequence of direct investigations of the occupation measures and their characterizations. For their arguments, extra conditions (e.g., bounded transition rate) are needed, and the role of the positive discount factor is important. In this connection, we point out that the discounted problem is a special case of the total undiscounted problem considered in the present paper, and the method of investigations here is quite different from \cite{DufourPiunovskiy:2016a}, and consequently, we do not need to impose any conditions on the growth of the transition and cost rates. More precisely, our investigation is based on the reduction of the gradual-impulsive CTMDP model to an equivalent but simpler standard CTMDP model. The reduction of gradual-impulsive model for piecewise deterministic processes to an equivalent model with gradual control only was proposed in \cite{DempsterYe:1995}. The reduction method in \cite{DempsterYe:1995} is different from the one here, and in fact, it induces a gradual control model with a much more complicated state space than the original one.

Our main contributions are as follows.
\begin{itemize}
\item[(a)] We fully justify that the gradual-impulsive CTMDP model can be reduced to an equivalent and simpler standard CTMDP model with the same state space. This reduction method was partially addressed and justified in \cite{PiunovskiyZhang:2021SICON}. The key difference is that in \cite{PiunovskiyZhang:2021SICON}, it was assumed that the transition intensities are strongly positive (separated from zero) at each state. This condition was essentially used in the argument in \cite{PiunovskiyZhang:2021SICON}. Here we manage to remove this extra condition, which, in our opinion, is a significant improvement. In fact, this turns out to be a delicate issue, and calls for a new and different proof. The new proof is based on the investigation of several new classes of control strategies, which can be of independent interest in their own right, and were not considered in \cite{PiunovskiyZhang:2021SICON}. The situation is much simpler if one only deals with strategies in simple form (e.g., stationary), but we consider general strategies.
\item[(b)] We establish the linear programming approach to solving constrained gradual-impulsive optimal control problem for CTMDPs with total undiscounted cost criteria. The linear program formulation itself is interesting, and was not reported in the previous literature, to the best of our knowledge. Moreover, no extra conditions on the growth of the transition and cost rates are needed. This is achieved by referring to the relevant results for the equivalent standard CTMDP problem, and thus also demonstrates the effectiveness of the reduction method fully justified in (a).
\end{itemize}

The rest of this paper is organized as follows. In Section \ref{PZ18Sec02} we describe the gradual-impulsive CTMDP model and the standard CTMDP model, and state the constrained optimal control problems under consideration. In Section \ref{2021June30Sec03} we present the main statements concerning the reduction method as well as the linear programming approach to the constrained optimal control problem. The justification of the reduction method is postponed to Section \ref{2021June30Sec05}, which also introduces some new classes of strategies and the auxiliary statements for them. The paper is ended with a conclusion in Section \ref{2021June30Sec06}. Some proofs are collected in the appendix.

\section{Model descriptions}\label{PZ18Sec02}
In this section, we describe the gradual-impulsive control model ${\cal M}$ and the model ${\cal M}^{GO}$ with gradual control only, as in \cite{PiunovskiyZhang:2021SICON}, which also goes back to \cite{Yushkevich:1980,Yushkevich:1988}.

We fist introduce some notations, definitions and facts to be used below, often without special reference. A Borel space is a Borel measurable subset of a complete separable metric space.  Suppose $\textbf{X}$ is a Borel space endowed with its Borel $\sigma$-algebra ${\cal B}(\textbf{X})$. Let ${\cal P}(\textbf{X})$ stand for the space of probability measures on $(\textbf{X}, {\cal B}(\textbf{X}))$.
We denote by ${\cal R}(\textbf{X})$ the collection of ${\cal P}(\textbf{X})$-valued measurable mappings on $(0,\infty)$ with any two elements therein being identified the same if they differ only on a null set with respect to the Lebesgue measure. Throughout this text, unless stated otherwise, by measurable we mean Borel measurable. For each $[-\infty,\infty]$-valued function $f$, $f^+$ and $f^-$ are its positive and negative parts. For brevity, by $f^{\pm}=g^{(\pm)}$ is meant $f^+=g^{(+)}$ and $f^-=g^{(-)}$ for some functions $g^{(+)}$ and $g^{(-)}$.

\subsection{Gradual-impulsive control model}\label{PZ18Subsection01}
We describe the primitives of the gradual-impulsive control model ${\cal M}$ as follows. The state space is $\textbf{X}$, the space of gradual controls is $\textbf{A}^G$, and the space of impulsive controls is $\textbf{A}^I$. It is assumed that $\textbf{X}$, $\textbf{A}^G$ and $\textbf{A}^I$ are all Borel spaces, endowed with their Borel $\sigma$-algebras ${\cal B}(\textbf{X}),$ ${\cal B}(\textbf{A}^G)$ and ${\cal B}(\textbf{A}^I)$, respectively. The transition rate, on which the gradual control acts, is given by $q(dy|x,a)$, which is a signed kernel from $\textbf{X}\times\textbf{A}^G$, endowed with its Borel $\sigma$-algebra, to ${\cal B}(\textbf{X}),$ satisfying the following conditions: $q(\Gamma|x,a)\in[0,\infty)$ for each $\Gamma\in{\cal B}(\textbf{X}), x\notin \Gamma;$
\begin{eqnarray*}
&&q(\textbf{X}|x,a)=0,~x\in \textbf{X},~a\in\textbf{A}^G;~\bar{q}_x:=\sup_{a\in \textbf{A}^G}q_x(a)<\infty,~x\in\textbf{X},
\end{eqnarray*}
where $q_x(a):=-q(\{x\}|x,a)$ for each $(x,a)\in\textbf{X}\times\textbf{A}^G.$ For notational convenience, we introduce
\begin{eqnarray*}
\tilde{q}(dy|x,a):=q(dy\setminus\{x\}|x,a),~\forall~x\in\textbf{X},~a\in\textbf{A}^G.
\end{eqnarray*}
If the current state is $x\in\textbf{X}$, and an impulsive control $b\in\textbf{A}^I$ is applied, then the state immediately following this impulse obeys the distribution given by $Q(dy|x,b)$, which is a stochastic kernel from $\textbf{X}\times\textbf{A}^I$ to ${\cal B}(\textbf{X}).$  We assume, without loss of generality that
\begin{eqnarray}\label{PZ18Assumption01}
Q(\{x\}|x,b)=0,~\forall~x\in\textbf{X},~b\in \textbf{A}^I.
\end{eqnarray} Finally, there are a family of cost rates and functions $\{c^G_i,c^I_i\}_{i=0}^J$, with $J$ being a fixed positive integer, representing the number of constraints in the concerned optimal control problem to be described below, see (\ref{PZ18Eqn18}). For each $i\in\{0,1,\dots,J\}$, $c^G_i$ and $c^I_i$ are $[-\infty,\infty]$-valued measurable functions on $\textbf{X}\times\textbf{A}^G$ and $\textbf{X}\times\textbf{A}^I$, respectively.


\begin{remark}
In what follows, we assume that $\textbf{A}^G$ and $\textbf{A}^I$ as two disjoint measurable subsets of a Borel space $\textbf{A}$ such that $\textbf{A}=\textbf{A}^G\cup\textbf{A}^I$. This is done without loss of generality, for otherwise, one can consider $\textbf{A}^G\times\{G\}$ instead of $\textbf{A}^G$ and $\textbf{A}^I\times\{I\}$ instead of $\textbf{A}^I$ and $ \textbf{A}=\textbf{A}^G\times\{G\}\cup \textbf{A}^I\times\{I\}$.
\end{remark}

The description of the system dynamics in the gradual-impulsive control problem is as follows. Assume $q_x(a)>0$ for each $x\in\textbf{X}$ and $a\in\textbf{A}^G$ for simplicity. At the initial time $0$ with the initial state $x_0$, the decision maker selects the triple $(\hat{c}_0,\hat{b}_0,\rho^0)$ with $\hat{c}_0\in[0,\infty]$, $\hat{b}_0\in\textbf{A}^I$, and {$\rho^0=\{\rho^0_t(da)\}_{t\in(0,\infty)} \in {\cal R}(\textbf{A}^G)$}. Then, the time until the next natural jump follows the nonstationary exponential distribution with the rate function $\int_{\textbf{A}^G} q_{x_0}(a)\rho^0_t(da)=:q_{x_0}(\rho^0_t)$. Here and below, unless stated otherwise, if $\rho\in {\cal R}(\textbf{A}^G)$, then $q_x(\rho_t):=\int_{\textbf{A}^G}q_x(a)\rho_t(da)$ and $\tilde{q}(dy|x,\rho_t):=\int_{\textbf{A}^G}\tilde{q}(dy|x,a)\rho_t(da).$ If by time $\hat{c}_0$, there is no occurrence of a natural jump, then the first sojourn time is $\hat{c}_0$, at which, the impulsive action $\hat{b}_0\in\textbf{A}^I$ is applied, and the next state $X_1$ follows the distribution $Q(dy|x_0,\hat{b}_0).$ If the first natural jump happens before $\hat{c}_0$, say at $t_1$, then the first sojourn time is $t_1$, and the next state $X_1$ follows the distribution
$
\frac{ \tilde{q}(dy|x_0, \rho^0_{t_1})}{ q_{x_0}( \rho^0_{t_1})}.
$
Except for the initial one, a decision epoch occurs immediately after a sojourn time. At the next decision epoch, the decision maker selects $(\hat{c}_1,\hat{b}_1,\rho^1)$, and so on. This leads to a natural description of the gradual-impulsive control problem as a discrete-time Markov decision process (DTMDP), which is presented next. This way of describing the gradual-impulsive control problem for a CTMDP goes back to Yushkevich \cite{Yushkevich:1988}.

The state space of the DTMDP model corresponding to the gradual-impulsive control problem
 is $\hat{\textbf{X}}:=\{(\infty,x_\infty)\}\cup [0,\infty)\times\textbf{X}$, where $(\infty,x_\infty)$ is an isolated point in $\hat{\textbf{X}}$. The first coordinate represents the previous sojourn time in the gradual-impulsive control problem, and the state of the controlled process in the gradual-impulsive control problem is given in the second coordinate. The inclusion of the first coordinate in the state allows us to consider control policies that select actions depending on the past sojourn times.

The action space of the DTMDP is $\hat{\textbf{A}}:=[0,\infty]\times \textbf{A}^I\times {\cal R}(\textbf{A}^G)$. Recall that ${\cal R}(\textbf{A}^G)$ is the collection of ${\cal P}(\textbf{A}^G)$-valued measurable mappings on $(0,\infty)$ with any two elements therein being identified the same if they differ only on a null set with respect to the Lebesgue measure, where ${\cal P}(\textbf{A}^G)$ stands for the space of probability measures on $(\textbf{A}^G, {\cal B}(\textbf{A}^G))$. We endow ${\cal P}(\textbf{A}^G)$ with its weak topology (generated by bounded continuous functions on $\textbf{A}^G$) and the Borel $\sigma$-algebra, so that ${\cal P}(\textbf{A}^G)$ is a Borel space, see Chapter 7 of \cite{Bertsekas:1978}. According to Lemma 3 of \cite{Yushkevich:1980}, each element in ${\cal R}(\textbf{A}^G)$ can be regarded as a stochastic kernel from $(0,\infty)$ to ${\cal B}(\textbf{A}^G)$. According to Lemma 1 of \cite{Yushkevich:1980}, the space ${\cal R}(\textbf{A}^G)$, endowed with the smallest $\sigma$-algebra with respect to which the mapping $\rho=(\rho_t(da))\in{\cal R}(\textbf{A}^G)\rightarrow \int_0^\infty e^{-t}g(t,\rho_t)dt$ is measurable for each bounded measurable function $g$ on $(0,\infty)\times {\cal P}(\textbf{A}^G)$, is a Borel space.

The transition probability $p$ in the DTMDP is defined as follows.
For each bounded measurable function $g$ on $\hat{\textbf{X}}$ and action $\hat{a}=(\hat{c},\hat{b},\rho)\in\hat{\textbf{A}}$,
\begin{eqnarray}\label{PZ18Eqn12}
&&\int_{\hat{\textbf{X}}}g(t,y)p(dt\times dy|(\theta,x),\hat{a})\nonumber\\
&:=&I\{\hat{c}=\infty\}\left\{ g(\infty,x_\infty)e^{-\int_0^\infty q_x(\rho_s)ds}+\int_{0}^\infty\int_{\textbf{X}}g(t,y)\tilde{q}(dy|x,\rho_t)e^{-\int_0^t q_x(\rho_s)ds} dt  \right\}\nonumber\\
&&+I\{\hat{c}<\infty\}\left\{\int_0^{\hat{c}} \int_{\textbf{X}} g(t,y)\tilde{q}(dy|x,\rho_t)e^{-\int_0^t q_x(\rho_s)ds}dt+e^{-\int_0^{\hat{c}} q_x(\rho_s)ds}\int_{\textbf{X}}g(\hat{c},y)Q(dy|x,\hat{b})  \right\}\nonumber\\
&=&\int_0^{\hat{c}} \int_{\textbf{X}} g(t,y)\tilde{q}(dy|x,\rho_t)e^{-\int_0^t q_x(\rho_s)ds}dt+I\{\hat{c}=\infty\} g(\infty,x_\infty)e^{-\int_0^\infty q_x(\rho_s)ds}\nonumber\\
&&+I\{\hat{c}<\infty\}e^{-\int_0^{\hat{c}} q_x(\rho_s)ds}\int_{\textbf{X}}g(\hat{c},y)Q(dy|x,\hat{b})
\end{eqnarray}
for each state
  $(\theta,x)\in[0,\infty)\times \textbf{X}$; and
$
\int_{\hat{\textbf{X}}}g(t,y)p(dt\times dy|(\infty,x_\infty),\hat{a}):=g(\infty,x_\infty).
$
The object $p$ defined above is indeed a stochastic kernel from $\hat{\textbf{X}}\times\hat{\textbf{A}}$ to ${\cal B}(\hat{\textbf{X}})$, see Lemma 2 of \cite{Yushkevich:1980} and its proof therein.
Similarly, the cost functions $\{l_i\}_{i=0}^J$ defined below are measurable on $\hat{\textbf{X}}\times\hat{\textbf{A}}\times\hat{\textbf{X}}$:
\begin{eqnarray}\label{PZ18Eqn11}
&&l_i((\theta,x),\hat{a},(t,y)):=I\{(\theta,x)\in [0,\infty)\times\textbf{X}\}\left\{\int_0^t c_i^G(x,\rho_s)ds+I\{t=\hat{c}<\infty\}c_i^I(x,\hat{b})\right\}\nonumber\\
&=&I\{x\in  \textbf{X}\}\left\{\int_0^t c_i^G(x,\rho_s)ds+I\{t=\hat{c}<\infty\}c_i^I(x,\hat{b})\right\}\nonumber\\
&:=&I\{x\in  \textbf{X}\}\left\{\int_0^t c_i^{G+}(x,\rho_s)ds+I\{t=\hat{c}<\infty\}c_i^{I+}(x,\hat{b})\right\}\nonumber\\
&&-I\{x\in  \textbf{X}\}\left\{\int_0^t c_i^{G-}(x,\rho_s)ds+I\{t=\hat{c}<\infty\}c_i^{I-}(x,\hat{b})\right\}\nonumber\\
&=:&l_i^{(+)}((\theta,x),\hat{a},(t,y))-l_i^{(-)}((\theta,x),\hat{a},(t,y)),
\end{eqnarray}
for each $i=0,1,\dots,J$ and $((\theta,x),\hat{a},(t,y))\in\hat{\textbf{X}}\times\hat{\textbf{A}}\times\hat{\textbf{X}}$. Here, $c_i^{G\pm},c_i^{I\pm}$ are the positive and negative parts of $c_i^G,c_i^I,$ and  the generic notation $\hat{a}=(\hat{c},\hat{b},\rho)\in\hat{\textbf{A}}$ of an action in this DTMDP model has been in use. The interpretation is that the pair $(\hat{c},\hat{b})$ is the pair of the planned time until the next impulse and the next planned impulse (provided that no natural jump occurs before then), and $\rho$ is (the rule of) the relaxed control to be used during the next sojourn time. Without loss of generality, the initial state is $(0,x_0)$, with some  $x_0\in \textbf{X}.$

Let $\{\hat{X}_n\}_{n=0}^\infty=\{(\hat{\Theta}_n,X_n)\}_{n=0}^\infty$ and $\{\hat{A}_n\}_{n=0}^\infty$ be the controlled and controlling process in this DTMDP model, and {$\{(\hat{C}_n,\hat{B}_n)\}_{n=0}^\infty$} the coordinate process corresponding to $\{(\hat{c}_n,\hat{b}_n)\}_{n=0}^\infty$ in $\{\hat{a}_n\}_{n=0}^\infty.$

Next, we define the concerned class of strategies in the gradual-impulsive control model.
\begin{definition}[(Ordinary) strategy in model ${\cal M}$]\label{PZ18Definition05}
Consider a sequence  $\sigma=\{\sigma_n\}_{n=0}^{\infty}$, where for each $n\ge 0$,
$\sigma_n$ is a stochastic kernel on ${\cal B}([0,\infty]\times \textbf{A}^I\times {\cal R}(\textbf{A}^G))$ given $
\hat{h}_n:=(\hat{x}_0,(\hat{c}_0,\hat{b}_0),\hat{x}_1,(\hat{c}_1,\hat{b}_1),\dots,\hat{x}_n).$
According to Proposition 7.27 of \cite{Bertsekas:1978} (or Proposition B.1.33 of \cite{PiunovskiyZhang:2020Book}),
\begin{eqnarray*}
\sigma_n(d\hat{c}\times d\hat{b}\times d\rho|\hat{h}_n)=\sigma_n^{(0)}(d\hat{c}\times d\hat{b}|\hat{h}_n)\sigma_n^{(1)}(d\rho|\hat{h}_n,\hat{c},\hat{b}),
\end{eqnarray*}
where $\sigma_n^{(0)}$ and $\sigma_n^{(1)}$ are some corresponding stochastic kernels.
If for each $n\ge 0$, there is a measurable mapping $\hat{F}_n$ mapping $(\hat{h}_n,\hat{c},\hat{b})$ to ${\cal R}(\textbf{A}^G)$ such that
\begin{eqnarray*}
\sigma_n^{(1)}(d\rho|\hat{h}_n,\hat{c},\hat{b})=\delta_{\hat{F}_n(\hat{h}_n,\hat{c},\hat{b})}(d\rho),
\end{eqnarray*}
then we call the sequence $\sigma=\{\sigma_n\}_{n=0}^\infty$, which is also identified with $\sigma=\{\sigma_n^{(0)},\hat{F}_n\}_{n=0}^\infty$, a strategy for the gradual-impulsive control model. The collection of all strategies for the gradual-impulsive CTMDP model is denoted by $\Sigma.$
\end{definition}

Note that the class of strategies defined above covers the particular case when one apriori determines a fixed time moment say $T$ of applying an impulse: this corresponds to $\sigma_n^{(0)}(d\hat{c}\times \textbf{A}^I|\hat{h}_n)=\delta_{T-\hat{t}_n}(d\hat{c})$ provided that $\hat{t}_n\le T,$ where $\hat{t}_n=\sum_{i=1}^n \hat{\theta}_n$ is the realized time of the $n$th jump moment, induced by either natural or active (impulsive) jumps.

\begin{definition}[Stationary strategy in model ${\cal M}$]\label{2021June30Definition05}
A strategy $\sigma=\{\sigma_n^{(0)},\hat{F}_n\}_{n=0}^\infty$ in model ${\cal M}$ is called stationary if for each $n\ge 0$,
\begin{eqnarray*}
 \sigma_n^{(0)}(d\hat{c}\times d\hat{b}|\hat{h}_n)=\sigma^{S,(0)}(d\hat{c}\times d\hat{b}|x_n), \hat{F}_n(\hat{h}_n,\hat{c},\hat{b})_t(da)=\hat{F}^S(x_n)(da),~\forall~t>0,
\end{eqnarray*}
where $\sigma^{S,(0)}(d\hat{c}\times d\hat{b}|x)$ and $\hat{F}^S(x)(da)$ are some stochastic kernels on ${\cal B}([0,\infty]\times \textbf{A}^I)$ concentrated on $\{0,\infty\}\times \textbf{A}^I$ and on ${\cal B}(\textbf{A}^G)$ given $x\in \textbf{X}$. We identify such a stationary strategy in ${\cal M}$ with $\sigma^S=(\sigma^{S,(0)},\hat{F}^S)$.
\end{definition}

Under a strategy $\sigma$ in the model ${\cal M}$, having in hand $\hat{h}_n$, the decision maker selects $(\hat{c}_n,\hat{b}_n)$ (possibly randomly), and after that, chooses
$\rho^n=\hat{F}_n(\hat{h}_n,\hat{c}_n,\hat{b}_n)$.

Given $\hat{x}_0=(0,x_0)\in\hat{\textbf{X}}$ and a strategy $\sigma$, let $\hat{ \rm P}^\sigma_{x_0}$ be the strategic measure in the DTMDP, and $\hat{{\rm E}}_{x_0}^\sigma$ the corresponding expectation. Then the concerned gradual-impulsive control problem with constraints reads
\begin{eqnarray}\label{PZ18Eqn18}
&&\mbox{Minimize over $\sigma\in\Sigma:$ } \hat{\rm E}_{x_0}^\sigma\left[\sum_{n=0}^\infty l_0(\hat{X}_n,\hat{A}_n,\hat{X}_{n+1})\right]=:\hat{W}_0(x_0,\sigma)\nonumber\\
&\mbox{subject to}& \hat{W}_j(x_0,\sigma):=\hat{\rm E}_{x_0}^\sigma\left[\sum_{n=0}^\infty l_j(\hat{X}_n,\hat{A}_n,\hat{X}_{n+1})\right]\le d_j,~j=1,\dots,J,
\end{eqnarray}
where $\{d_j\}_{j=1}^J\subset\mathbb{R}^{J}$ is a fixed vector of constants, $x_0$ is a fixed element of $\textbf{X}$, and
\begin{eqnarray*}
\hat{\rm E}_{x_0}^\sigma\left[\sum_{n=0}^\infty l_i(\hat{X}_n,\hat{A}_n,{\hat{X}_{n+1}})\right]:=\hat{\rm E}_{x_0}^\sigma\left[\sum_{n=0}^\infty l_i^{(+)}(\hat{X}_n,\hat{A}_n,\hat{X}_{n+1})\right]-\hat{\rm E}_{x_0}^\sigma\left[\sum_{n=0}^\infty l_i^{(-)}(\hat{X}_n,\hat{A}_n,\hat{X}_{n+1})\right]
\end{eqnarray*}
with $\infty-\infty:=\infty$ being adopted here, and $l_i^{(\pm)}$ being defined in (\ref{PZ18Eqn11}).


\subsection{Standard CTMDP model}
In a standard CTMDP model, there is only gradual control, which is selected according to relaxed policies\footnote{The term policy is a synonym of the term strategy, but we use ``policy'' exclusively for models with gradual control only.} . Its system primitives are the following objects
\begin{eqnarray*}
{\cal M}^{GO}:=\{\textbf{X},\textbf{A},q^{GO},\{c_i^{GO}\}_{i=0}^J\}.
\end{eqnarray*}
Here the state and action spaces $\textbf{X}$ and $\textbf{A}$ are Borel spaces,   $q^{GO}$ is the transition rate from $\textbf{X}\times\textbf{A}$ to ${\cal B}(\textbf{X})$, and $\{c_i^{GO}\}_{i=0}^J$ is the collection of measurable functions on $\textbf{X}\times\textbf{A}$, representing the cost rates, $J\ge 0$ is a fixed integer.  The superscript ``$GO$'' abbreviates ``gradual only'', as the model only allows gradual controls.

In the standard CTMDP model ${\cal M}^{GO} $, a decision epoch occurs after each natural jump of the controlled process (except for the initial decision epoch at time zero). At each decision epoch, one selects the relaxed control function $\rho\in\cal R(\textbf{A})$ until the next decision epoch occurs. We sketch the more rigorous construction as follows.
The  sample space $\Omega$ is taken as the union of $ (\textbf{X}\times   (0,\infty))^\infty$ and the collection of sequences in the form $( x_0, \theta_1,x_1, \dots,\theta_{m-1},x_{m-1}, \infty, x_\infty, \infty,x_\infty, \dots)$, where $m\ge 1$, and $x_\infty\notin\textbf{X}$ is an isolated point. We endow $\Omega$ with the $\sigma$-algebra ${\cal F}$ obtained as the trace of ${\cal B}(  (\textbf{X}_\infty\times (0,\infty])^\infty)$ on $\Omega$, where $\textbf{X}_\infty=\textbf{X}\cup\{x_\infty\}$.
The generic notation for an element of $\Omega$ is $\omega.$ For each $\omega\in\Omega$, define $\theta_0:=0$, $t_n:=\sum_{i=0}^{n}\theta_i,$ $h_n:=(x_0,\theta_1,x_1,\dots,\theta_n,x_n)$ for each $n\ge 0.$ The collection of all possible $h_n$ is denoted as $\textbf{H}_n$ for each $n\ge 0.$ Let us put $t_\infty:=\lim_{n\rightarrow \infty} t_n$, which exists.
When regarded as coordinate variables, we use capital letters $\Theta_n,$ $T_n$, $X_n,$ and $H_n$ corresponding to $\theta_n,t_n,x_n$ and $h_n$. The state process $\{X(t)\}_{t\ge 0}$ is defined by
$
X(t):=X_n
$ if $T_n\le t<T_{n+1}$ for some $n\ge 0,$ and $X(t):=x_\infty$ if $t\ge T_\infty.$
As usual, we omit $\omega$ whenever the context excludes confusion.

\begin{definition}[(Ordinary) policy in ${\cal M}^{GO}$]
A policy $\overline{S}$ in the standard CTMDP model ${\cal M}^{GO}$ is the following object: $\overline{S}= \{\overline{F}_n\}_{n=0}^\infty $,  for each $n\ge 0,$ $\overline{F}_n$ is a measurable mapping on $\textbf{H}_n$ taking values in ${\cal R}(\textbf{A})$.
\end{definition}

\begin{definition}[Markov policy in ${\cal M}^{GO}$]\label{2021June30Definition01}
A policy $\overline{S}=\{\overline{F}_n\}_{n\ge 0}$ in ${\cal M}^{GO}$ is called Markov if $\overline{F}_n(h_n)=\overline{F}_n^M(x_n)$ for some measurable mapping $\overline{F}_n^M$ from $\textbf{X}$ to ${\cal R}(\textbf{A})$. In this case, we identify $\overline{S}$ with $\{\overline{F}_n^M\}_{n\ge 0}=:\overline{S}^M.$
\end{definition}

\begin{definition}[Stationary policy in model ${\cal M}^{GO}$]\label{2021June30Definition02}
A policy $\overline{S}=\{\overline{F}_n\}_{n\ge 0}$ in ${\cal M}^{GO}$ is called stationary if $\overline{F}_n(h_n)_t(da)=\overline{F}^S(x_n)(da)$ for some stochastic kernel $\overline{F}^S(x)(da)$ on ${\cal B}(\textbf{A})$ given $x\in \textbf{X}$. In this case, we identify such a stationary policy $\overline{S}$ with $\overline{F}^S.$
\end{definition}

\begin{remark}\label{PZ18Remark01}
We put $q^{GO}_{x_\infty}(a)\equiv 0\equiv q^{GO}(\Gamma|x_\infty,a)$ for all $\Gamma\in{\cal B}(\textbf{X})$ and $c^{GO}_i(x_\infty,a)\equiv 0$ in what follows.
\end{remark}

Given a policy $\overline{S}= \{\overline{F}_n\}_{n=0}^\infty $ and initial state $x_0\in\textbf{X}$, there is a unique probability measure ${\rm P}_{x_0}^{\overline{S}}$ on $(\Omega,{\cal F})$ such that ${\rm P}_{x_0}^{\overline{S}}(X_0\in dx)=\delta_{x_0}(dx)$, and for each $n\ge 1$ and  $\Gamma_1\in {\cal B}([0,\infty))$, $\Gamma_2\in {\cal B}(\textbf{X})$,
\begin{eqnarray*}
&&{\rm P}_{x_0}^{\overline{S}}(  \Theta_n\in \Gamma_1,~X_n\in \Gamma_2|H_{n-1})\\
&=& \int_{\Gamma_1}e^{-\int_0^s q^{GO}_{X_{n-1}}(\overline{F}_{n-1}(H_{n-1})_t)dt}\tilde{q}^{GO}(\Gamma_2|{X_{n-1},\overline{F}_{n-1}(H_{n-1})_s)}ds ;\\
&&{\rm P}_{x_0}^{\overline{S}}(\Theta_n=\infty,~X_n=x_\infty|H_{n-1})=e^{-\int_0^\infty q^{GO}_{X_{n-1}}(\overline{F}_{n-1}(H_{n-1})_t)dt};
\end{eqnarray*}
and
\begin{eqnarray*}
{\rm P}_{x_0}^{\overline{S}}(\Theta_n=\infty,~X_n\in \Gamma_2|H_{n-1})= {\rm P}_{x_0}^{\overline{S}}(\Theta_n\in\Gamma_1,~X_n=x_\infty|H_{n-1})=0.
\end{eqnarray*}
Let the expectation corresponding to ${\rm P}^{\overline{S}}_{x_0}$ be denoted as ${\rm E}_{x_0}^{\overline{S}}$. We consider the following optimal control problem corresponding to problem (\ref{PZ18Eqn18}):
\begin{eqnarray}\label{PZ18Eqn19Bar}
&&\mbox{Minimize over $\overline{S}:$ }  W_0(x_0,\overline{S}):= {\rm E}_{x_0}^{\overline{S}}\left[\sum_{n=0}^\infty I\{T_n<\infty\}\int_{T_n}^{T_{n+1}} c^{GO}_0({X}_n,\overline{F}_{n}(H_{n})_{t-T_n})dt\right]  \nonumber\\
&\mbox{subject to}&   W_j(x_0,\overline{S}):={\rm E}_{x_0}^{\overline{S}}\left[\sum_{n=0}^\infty  I\{T_n<\infty\} \int_{T_n}^{T_{n+1}} c^{GO}_j({X}_n,\overline{F}_{n}(H_{n})_{t-T_n})dt\right] \le d_j,\nonumber\\
&&~j=1,\dots,J,
\end{eqnarray}
where \begin{eqnarray*}
 &&{\rm E}_{x_0}^{\overline{S}}\left[I\{T_n<\infty\}\sum_{n=0}^\infty \int_{T_n}^{T_{n+1}} c^{GO}_i({X}_n,\overline{F}_{n}(H_{n})_{t-T_n}) dt\right] \\
&:=&
{\rm E}_{x_0}^{\overline{S}}\left[\sum_{n=0}^\infty I\{T_n<\infty\}\int_{T_n}^{T_{n+1}} {c^{GO}}^+_i({X}_n,\overline{F}_{n}(H_{n})_{t-T_n}) dt\right]\\
 &&-{\rm E}_{x_0}^{\overline{S}}\left[\sum_{n=0}^\infty I\{T_n<\infty\} \int_{T_n}^{T_{n+1}} {c^{GO}}^-_i({X}_n,\overline{F}_{n}(H_{n})_{t-T_n})dt\right],
\end{eqnarray*}
with $\infty-\infty:=\infty$ being accepted and $c^{GO\pm}_i$ being the positive and negative part of $c_i^{GO}$, respectively.  Here, the constants $J$ and $\{d_j\}_{j=1}^J$ are the same as in problem (\ref{PZ18Eqn18}), and we have used the following notation: for each probability measure $\mu$ on ${\cal B}(\textbf{X})$ and measurable function $f$ on $\textbf{X}$, we put $f(\mu):=\int_{\textbf{X}}f(x)\mu(dx)$ whenever the right hand side is well defined. This notation is only for brevity, and will be used when there is no potential confusion regarding the underlying space $\textbf{X}$.

For the future, it is convenient to note that we may also write
\begin{eqnarray*}
W_i(x_0,\overline{S})&=&{\rm E}_{x_0}^{\overline{S}}\left[\sum_{n=0}^\infty  I\{X_n\in \textbf{X}\} \int_{0}^{\Theta_{n+1}} {c^{GO}}^+_i({X}_n,\overline{F}_{n}(H_{n})_{t}) dt\right]\\
&&-{\rm E}_{x_0}^{\overline{S}}\left[\sum_{n=0}^\infty  I\{X_n\in \textbf{X}\} \int_{0}^{\Theta_{n+1}} {c^{GO}}^-_i({X}_n,\overline{F}_{n}(H_{n})_{t})dt\right].
\end{eqnarray*}

\section{Main results}\label{2021June30Sec03}
\subsection{Reduction results}

In the rest of this paper, we consider the following standard CTMDP model ${\cal M}^{GO}$ induced by the gradual-impulsive control model ${\cal M}$, defined as follows
\begin{eqnarray}\label{PZ18Eqn10}
&&\textbf{A}:=\textbf{A}^I\cup\textbf{A}^G;~q^{GO}(dy|x,a):=q(dy|x,a),~\forall~(x,a)\in\textbf{X}\times\textbf{A}^G;\nonumber\\
&&\tilde{q}^{GO}(dy|x,a):=Q(dy|x,a),~q^{GO}_x(a):=1,~\forall~(x,a)\in\textbf{X}\times\textbf{A}^I;\nonumber\\
&& c^{GO}_i(x,a):=c_i^G(x,a),~\forall~(x,a)\in\textbf{X}\times\textbf{A}^G;~c^{GO}_i(x,a):=c_i^I(x,a),~\forall~(x,a)\in\textbf{X}\times\textbf{A}^I.\nonumber\\
\end{eqnarray}
(Equality (\ref{PZ18Assumption01}) guarantees that $q^{GO}$ defined in the above is indeed a transition rate.)

\begin{definition}
A policy (or strategy) in a model is said to replicate another policy (or strategy) in a possibly different model if the system performances of the two policies or strategies in their respective models coincide.
\end{definition}
We say that the gradual-impulsive control model ${\cal M}$ can be reduced to the model ${\cal M}^{GO}$ with gradual control if each strategy in ${\cal M}$ is replicated by a policy in ${\cal M}^{GO}$, and each policy in ${\cal M}^{GO}$ is replicated by a strategy in ${\cal M}.$

One purpose of this section is to show that the gradual-impulsive control model ${\cal M}$ can be reduced to the model ${\cal M}^{GO}$ with gradual control only.

\begin{theorem}\label{2021June30Theorem01}
The gradual-impulsive control model ${\cal M}$ can be reduced to the model ${\cal M}^{GO}$ with gradual control only. That is, each policy in the gradual control model ${\cal M}^{GO}$ can be replicated by a strategy in the gradual-impulsive control model ${\cal M}^{GO},$ and vice versa.
\end{theorem}
The proof of this theorem is postponed to Section \ref{2021June30Sec05}.

Here let us provide some comments and discussions. Such a reduction result is desirable and useful, because, on the one hand, the gradual-impulsive control model considered here is rather general (in particular, the impulse can be applied at any time moment determined a priori, and the optimal stopping problem is a special case), and investigations of such models following a direct method can be involving (processes are not stochastically continuous, multiple impulses are allowed at a single time leading to rather complicated states, etc), see \cite{DufourPiunovskiy:2016,DufourPiunovskiy:2016a,Yushkevich:1988}; on the other hand, the theory for standard CTMDP models is fairly matured: for some recent monographs, see, e.g., \cite{Guo:2009,PiunovskiyZhang:2020Book,Prieto-RumeauHernandez-Lerma:2012}.

This reduction issue was partially addressed in \cite{PiunovskiyZhang:2021SICON}. Indeed, it was established in  Theorem 3.2 of \cite{PiunovskiyZhang:2021SICON} that any strategy in ${\cal M}$ can be replicated by a policy in ${\cal M}^{GO}.$
The opposite direction is more delicate. The corresponding statement, collected as Proposition \ref{PZ18Theorem01} below, was established in \cite{PiunovskiyZhang:2021SICON} under the following extra condition:
\begin{condition}\label{2021June30Con01}
For each $x\in\textbf{X},$ there is some $\epsilon>0$ such that $q_x(a)\ge \epsilon>0$ for all $a\in\textbf{A}^G.$
\end{condition}


\begin{proposition}\label{PZ18Theorem01}
Suppose Condition \ref{2021June30Con01} is satisfied. Then each policy in the gradual control model ${\cal M}^{GO}$ can be replicated by a strategy in the gradual-impulsive control model ${\cal M}^{GO},$ i.e., for each policy $\overline{S}$ in the gradual control model ${\cal M}^{GO}$ there is a strategy $\sigma$ in the gradual-impulsive control model ${\cal M}$ such that $\hat{W}_i(x_0,\sigma)=W_i(x_0,\overline{S}).$
\end{proposition}
\par\noindent\textit{Proof.} See Theorem 3.1 of \cite{PiunovskiyZhang:2021SICON}. $\hfill\Box$ \bigskip

A main contribution of this paper lies in showing that Condition \ref{2021June30Con01} can be withdrawn from Proposition \ref{PZ18Theorem01}, and that removal would also complete the proof of Theorem  \ref{2021June30Theorem01}. We underline that the argument in the proof of Theorem 3.1 of \cite{PiunovskiyZhang:2021SICON} essentially made use of Condition \ref{2021June30Con01}. Here we will develop a different method, based on investigations of auxiliary (new) classes of control strategies and policies for the model ${\cal M}$ and for the model ${\cal M}^{GO}$, which are introduced in Section \ref{2021June30Sec05}, where relevant properties of the introduced auxiliary classes of strategies are presented and can be of independent interest. They were not considered in \cite{PiunovskiyZhang:2021SICON}.

The situation is simpler if we consider stationary policies in model ${\cal M}^{GO}$. They can be indeed replicated by stationary strategies in model ${\cal M}$ without Condition \ref{2021June30Con01}, as observed in the next statement. Its proof can be done directly without involving auxiliary strategies, though the argument cannot handle the case of general strategies.
\begin{proposition}\label{2021June30Proposition01}
Each stationary policy  $\overline{F}^S$ in ${\cal M}^{GO}$ is replicated by the stationary strategy $\sigma^S=(\sigma^{S,(0)},\hat{F}^S)$  defined as follows: for each $x\in O$ with
\begin{eqnarray*}
O:=\left\{x\in \textbf{X}:~ \int_{\textbf{A}^G} q_x(a) \overline{F}^S(x)(da) +  \overline{F}^S(x)(\textbf{A}^I)>0\right\},
\end{eqnarray*}
On ${\cal B}(\textbf{A}^I)$:
\begin{eqnarray*}
&&\sigma^{S,(0)}(\Gamma\times d\hat{b}|x)=0~\forall~\Gamma\in{\cal B}(0,\infty),\\
&&\sigma^{S,(0)}(\{0\}\times d\hat{b}|x)= \frac{\overline{F}^S(x)(d\hat{b})}{\int_{\textbf{A}^G} q_x(a) \overline{F}^S(x)(da) +  \overline{F}^S(x)(\textbf{A}^I)},\\
&&\sigma^{S,(0)}(\{\infty\}\times d\hat{b}|x)=
p^{\ast\ast}(d\hat{b}) \frac{\int_{\textbf{A}^G} q_x(a) \overline{F}^S(x)(da) }{\int_{\textbf{A}^G} q_x(a) \overline{F}^S(x)(da) +  \overline{F}^S(x)(\textbf{A}^I)}
\end{eqnarray*}
where $p^{\ast\ast}\in {\cal P}(\textbf{A}^I)$ is an arbitrarily fixed probability measure on ${\cal B}(\textbf{A}^I)$;
\begin{eqnarray*}
\hat{F}^S(x)(da)=
\begin{cases}
\frac{\overline{F}^S(x)(da\cap \textbf{A}^G)}{\overline{F}^S(x)(\textbf{A}^G)}&\mbox{~if $\overline{F}^S(x)(\textbf{A}^G)>0$,}\\
p^\ast(da)&\mbox{~otherwise,}
\end{cases}
\end{eqnarray*}
where $p^{\ast}\in {\cal P}(\textbf{A}^G)$ is an arbitrarily fixed probability measure on ${\cal B}(\textbf{A}^G);$ whereas for each $x\in \textbf{X}\setminus O,$
\begin{eqnarray*}
\sigma^{S,(0)}(d\hat{c}\times d\hat{b}|x)=\delta_{\infty}(d\hat{c})p^{\ast\ast}(d\hat{b}),~
\hat{F}^S(x)(da)=
\overline{F}^S(x)(da).
\end{eqnarray*}
\end{proposition}
The proof of Proposition \ref{2021June30Proposition01} is given in the appendix.

\subsection{Optimality results}


In this subsection, we firstly impose a compactness-continuity condition. Then, under that condition, we may conclude the existence of an optimal stationary strategy (out of the class of ordinary strategies) for problem (\ref{PZ18Eqn18}), and establish a linear program, solving which, one can produce the optimal stationary strategy. This is achieved by making use of known results for the standard CTMDP problem (\ref{PZ18Eqn19Bar}) together with the reduction results in the previous subsection. For this reason, in its proof, we primarily refer the reader to the corresponding references for standard CTMDPs instead of full details.  This linear program approach for problem (\ref{PZ18Eqn18}) is in its own right of interest. In the current general form, it was not reported in the literature, to the best of our knowledge.

\begin{condition}\label{PZ18Condition02}
\begin{itemize}
\item[(a)] $\textbf{A}^G$ and $\textbf{A}^I$ are compact.
\item[(b)]
The functions $\{c_i^G\}_{i=0}^J$ and $\{c_i^I\}_{i=0}^J$ are $[0,\infty]$-valued and lower semicontinuous on $\textbf{X}\times\textbf{A}^G$ and $\textbf{X}\times\textbf{A}^I$, respectively.
\item[(c)] For each bounded continuous function $f$ on $\textbf{X}$, the functions $(x,a)\in \textbf{X}\times\textbf{A}^G\rightarrow \int_{\textbf{X}}f(y)\tilde{q}(dy|x,a)$ and $(x,b)\in \textbf{X}\times\textbf{A}^I\rightarrow \int_{\textbf{X}}f(y)Q(dy|x,b)$ are continuous.
    \end{itemize}
\end{condition}

Under Condition \ref{PZ18Condition02}, we present the linear program formulation, for which some additional notations are introduced.
Let $v^\ast$ be the minimal nonnegative lower semicontinuous function on $\textbf{X}$ satisfying the first equality in
\begin{eqnarray*}
v^\ast(x)&=&\inf_{a\in \textbf{A}}\left\{\frac{\sum_{j=0}^J c_j^{GO}(x,a)}{\epsilon+q^{GO}_x(a)}+\frac{\int_{\textbf{X}}v^\ast(y) \tilde{q}^{GO}(dy|x,a) +\epsilon v^\ast(x)}{\epsilon+q^{GO}_x(a)}\right\}\\
&=&\frac{\sum_{j=0}^J c_j^{GO}(x,f^\ast(x))}{\epsilon+q^{GO}_x(f^\ast(x))}+\frac{\int_{\textbf{X}}v^\ast(y) \tilde{q}^{GO}(dy|x,f^\ast(x))+\epsilon v^\ast(x)}{\epsilon+q^{GO}_x(f^\ast(x))} ,~x\in\textbf{X}
\end{eqnarray*}
(recall $\textbf{A}=\textbf{A}^G\cup\textbf{A}^I$), where $f^\ast$ is a measurable mapping from $\textbf{X}$ to $\textbf{A}.$ Note that $v^\ast$ is actually independent of $\epsilon>0$, and the existence of $v^\ast$ and $f^\ast$ is guaranteed under Condition \ref{PZ18Condition02}, according to, e.g., Theorem 4.2.1 of \cite{PiunovskiyZhang:2020Book} and its proof.  Put $\textbf{R}:=\{x\in \textbf{X}:~ v^\ast(x)>0\}$. (The intuitive meaning of $\textbf{R}^c$ is the part of the state space, at which it is optimal to apply $f^\ast$ in the model ${\cal M}^{GO}$: the process will remain there with no cost being incurred. Thus, the nontrivial part is to determine the control in ${\cal M}^{GO}$when the process is in $\textbf{R}$.)  Then consider the following linear program:
\begin{eqnarray}\label{2021June30Eqn38}
 &&\int_{\textbf{R}\times \textbf{A}^G} c_0^G(x,a)\nu(dx\times da)+\int_{\textbf{R}\times \textbf{A}^I} c_0^I(x,a)\nu(dx\times da)\rightarrow \min_{\nu}\nonumber\\
&\mbox{s.t.}&~  \int_{\textbf{A}^G}q_y(a)\nu(dy\times da)+  \nu(dx\times \textbf{A}^I)= \delta_{x_0}(dx)+ \int_{\textbf{R}\times\textbf{A}^G}\tilde{q}(dx|y,a)\nu(dy\times da)\nonumber\\
&&~~~~~~~~~~+\int_{\textbf{R}\times\textbf{A}^I}Q(dx|y,a)\nu(dy\times da);  \\
&&\int_{\textbf{R}\times \textbf{A}^G} c_j^G(x,a)\nu(dx\times da)+\int_{\textbf{R}\times \textbf{A}^I} c_j^I(x,a)\nu(dx\times da)\le d_j,~j\in\{1,2,\dots,J\};\nonumber\\
&&\mbox{$\nu$ is a measure on ${\cal B}(\textbf{R}\times\textbf{A}):$}\nonumber\\
&&\mbox{$\nu(dx\times\textbf{A})$ is a $\sigma$-finite measure on ${\cal B}(\textbf{R});$}\nonumber\\
&&\mbox{~ $\int_{\textbf{A}^G}q_x(a)\nu(dx\times da)+ \nu(dx\times \textbf{A}^I)$ is $\sigma$-finite on ${\cal B}(\textbf{R}).$}\nonumber
\end{eqnarray}

\begin{theorem}[Linear programming approach]\label{2021June30Theorem06}
Suppose that Condition \ref{PZ18Condition02} is satisfied, and there is a feasible strategy $\sigma$ for problem (\ref{PZ18Eqn18}) such that $\hat{W}_0(x_0,\sigma)<\infty.$ Then the following assertions hold.
\begin{itemize}
\item[(a)] There exists an optimal stationary strategy for problem (\ref{PZ18Eqn18}).
\item[(b)]  If the linear program (\ref{2021June30Eqn38}) has a feasible solution, which is the case if problem (\ref{PZ18Eqn18}) has a feasible strategy with finite value, then,  the linear program has an optimal solution, say $\nu^\ast$. Consider the stochastic kernel $\overline{F}^S(x)(da)$ on ${\cal B}(\textbf{A})$ given $x\in\textbf{X}$ satisfying $\nu^\ast(\Gamma\times da)=\nu^\ast(dx\times \textbf{A})\overline{F}^S(x)(da)$ for each $\Gamma\in {\cal B}(\textbf{R})$, and $\overline{F}^S(x)(da)=\delta_{f^\ast(x)}(da)$ for each $x\in \textbf{X}\setminus\textbf{R}.$ (Such a stochastic kernel exists because $\nu^\ast(dx\times \textbf{A})$ is $\sigma$-finite on ${\cal B}(\textbf{R})$.) Then the stationary strategy $\sigma^S=(\sigma^{S,(0)},\hat{F}^S)$ defined in terms of $\overline{F}^S$ in Proposition \ref{2021June30Proposition01} is optimal for problem (\ref{PZ18Eqn18}).
\end{itemize}
\end{theorem}

\par\noindent\textit{Proof.}
(a) By Theorem \ref{2021June30Theorem01}, the gradual-impulsive optimal control problem (\ref{PZ18Eqn18}) can be reduced to the standard CTMDP problem (\ref{PZ18Eqn19Bar}) with gradual control only. Statement (a) follows from this reduction, Theorem 4.2.2(b) of \cite{PiunovskiyZhang:2020Book}, and Proposition \ref{2021June30Proposition01}.

(b) The induced standard CTMDP problem (\ref{PZ18Eqn19Bar}) can be reduced to a discrete-time Markov decision process (DTMDP) (without any compactness-continuity conditions). For the details, see Theorems 4.2.1 and 6.2.1 of \cite{PiunovskiyZhang:2020Book}. Now statement (b) follows from this chain of reductions, Proposition \ref{2021June30Proposition01}, and the relevant result for the induced DTMDP problem obtained in \cite{Dufour:2012}, which is also collected in Proposition C.2.18 and Remark C.2.4 of \cite{PiunovskiyZhang:2020Book}. See the proof of Theorem 4.2.2 of \cite{PiunovskiyZhang:2020Book} for the details. $\hfill\Box$
\bigskip

The linear programming approach for problem (\ref{PZ18Eqn18}) was not reported in the previous literature. Theorem \ref{2021June30Theorem06} can be viewed as a significant extension of the corresponding result in \cite{DufourPiunovskiy:2016a}. Only the discounted model was considered in \cite{DufourPiunovskiy:2016a}, which, as we underline, follows a different method to obtain the linear program formulation. The method in \cite{DufourPiunovskiy:2016a}, on the one hand, requires extra conditions on the growth (boundedness) of the transition rate, makes use the presence of discounting, and does not involve any reduction to standard CTMDP models, on the other hand. In this connection, we point out that the discounted problem is a special case of the total undiscounted problem (\ref{PZ18Eqn18}) considered here: the justification can be found in Subsection 7.3.1 (in particular, Theorem 7.3.1) of \cite{PiunovskiyZhang:2020Book}. Specialized to discounted problems, the linear program (\ref{2021June30Eqn38}) is consistent with the linear program established in \cite{DufourPiunovskiy:2016a} (see, in particular, the equalities in the proof of Theorem 4.6 therein).

\begin{definition}[Deterministic stationary strategy in model ${\cal M}$]
A stationary strategy $\sigma^S=(\sigma^{S,(0)},\hat{F}^S)$ in model ${\cal M}$ is called deterministic stationary if
\begin{eqnarray*}
 \sigma^{S,(0)}(d\hat{c}\times d\hat{b}|x)=\delta_{\varphi(x)}(d\hat{c})\delta_{\zeta(x)}(d\hat{b}),~ \hat{F}^S(x)(da)=\delta_{f^S(x)}(da),~\forall~t>0,
\end{eqnarray*}
where  $\varphi$ (or $\zeta$, $f^S$) is a measurable mapping from $\textbf{X}$ to $\{0,\infty\}$ ($\textbf{A}^I$, $\textbf{A}^G$, respectively). We identify such a deterministic stationary strategy in ${\cal M}$ with $(\varphi,\zeta,f^S)$.
\end{definition}

The next example demonstrates that deterministic stationary strategies are not sufficient  for the constrained problem (\ref{PZ18Eqn18}).
\begin{example}
Let $\textbf{X}=\{0,1,2,\dots\}$, $\textbf{A}^G=\{a\}$, $\textbf{A}^I=\{b\}$ with $a\ne b$, so that we may put $\textbf{A}=\{a,b\}=\textbf{A}^G\cup\textbf{A}^I$. Let $q_0(a)=1=q(\{1\}|0,a)$, $q_x(a)=0$ for all $x\in\{1,2,\dots\}$, $Q(\{x+1\}|x,b)=1$ for all $x\in \textbf{X}$.  Finally, fix $J=1$, $d_1=1$, $x_0=0,$ and consider the cost rates and functions defined by
\begin{eqnarray*}
&&c_0^G(0,a)=1,~c_0^G(x,a)=0~\forall~x\in\{1,2,\dots\}; \\
&&c_0^I(x,b)=0~\forall~x\in\{0,1,2,\dots\};\\
&&c_1^G(x,a)=0~\forall~x\in\{0,1,2,\dots\};\\
&&c_1^I(0,b)=2,~c_1^I(x,b)=0~\forall~x\in\{1,2,\dots\}.
\end{eqnarray*}
Apparently, since the process is essentially only controlled at the state $x=0,$ (once the process leaves the state $0$, no further cost will be incurred), as far as the performance of deterministic stationary strategies is concerned, one only needs to consider deterministic  stationary strategies in the following form: $\sigma^{DS}=(\varphi,\zeta,f^S)$ given by $\varphi(0)=0$ and $\sigma^{'DS}=(\varphi',\zeta,f^S)=\varphi'(0)=\infty.$
We may compute
\begin{eqnarray*}
&&\hat{W}_0(0,\sigma^{DS})=0,~\hat{W}_1(0,\sigma^{DS})=2>d_1=1;\\
&&\hat{W}_0(0,\sigma^{'DS})=1,~\hat{W}_1(0,\sigma^{'DS})=0.
\end{eqnarray*}
Consequently, $\sigma^{DS}$ is not feasible for problem (\ref{PZ18Eqn18}). Now consider $\sigma^S=(\sigma^{S,(0)},\hat{F}^S)$ such that $\sigma^{S,(0)}(\{0\}\times\{b\}|0)=0.5=\sigma^{S,(0)}(\{\infty\}\times\{b\}|0)$. Then one can verify that
\begin{eqnarray*}
\hat{W}_0(0,\sigma^{S})=\frac{1}{2}<\hat{W}_0(0,\sigma^{'DS}),~\hat{W}_1(0,\sigma^{S})=1,
\end{eqnarray*}
which is feasible and strictly outperforms $\sigma^{'DS}$, and thus strictly outperforms any feasible deterministic stationary strategy.
\end{example}

\section{Auxiliary statements and proof of Theorem \ref{2021June30Theorem01}}\label{2021June30Sec05}

The proof of Theorem \ref{2021June30Theorem01} goes in several steps, and, as was aforementioned, makes use of auxiliary classes of strategies in the model ${\cal M}$ and policies in the model ${\cal M}^{GO}$, which are introduced in separate subsections.

\subsection{Pseudo-Poisson-related policy in the model ${\cal M}^{GO}$ with gradual control only}
In what follows, we fix some strictly positive constant $\lambda>0.$

Let \begin{eqnarray*}
{\bf \Xi}^{GO}:=[0,\infty)\times \textbf{A}\times ((0,\infty]\times \textbf{A})^\infty
\end{eqnarray*} be the countable product.
The generic notation for an element of ${\bf \Xi}^{GO}$ is $\xi=\{(\psi_n,\alpha_n)\}_{n\ge 0}\in {\bf \Xi}^{GO}$.  Consider the coordinate random variables (viewing $({\bf \Xi}^{GO},{\cal B}({\bf \Xi}^{GO}))$ as a sample space): for each $\xi=\{(\psi_n,\alpha_n)\}_{n\ge 0}\in {\bf \Xi}^{GO}$, $\Psi_n(\xi):=\psi_n$ and $\Phi_n(\xi):=\alpha_n$, and $\tau_{n}:=\sum_{k=0}^n \psi_k.$

Let \begin{eqnarray*}
\overline{\lambda}(a):=\lambda I\{a\in \textbf{A}^G\}~\forall~a\in\textbf{A}.
\end{eqnarray*}

\begin{definition}[Pseudo-Poisson-related policy in ${\cal M}^{GO}$]
A pseudo-Poisson-related policy in ${\cal M}^{GO}$ is given by a sequence of stochastic kernels $\overline{S}^P=\{\overline{p}_n(d\xi|x)\}_{n\ge 0}$ on ${\cal B}({\bf \Xi}^{GO})$ from $x\in \textbf{X}$, where for each $n\ge 0$ and $x\in\textbf{X}$, under $\overline{p}_n(d\xi|x)$,
\begin{eqnarray*}
&&\overline{p}_n(\Psi_0\in dt|x)=\delta_0(dt),
\end{eqnarray*}
and the random vectors $(\Phi_0,\Psi_1),(\Phi_1,\Psi_2)),\dots$ are mutually independent satisfying
\begin{eqnarray*}
&& \overline{p}_n(\Phi_k\in da|x)=:\overline{p}_{n,k}(da|x)~\forall~k\in\{0, 1,2,\dots\},\\
&&\overline{p}_n(\Phi_k\in da,~\Psi_{k+1}> t|x)=e^{-\overline{\lambda}(a) t}\overline{p}_{n,k}(da|x),~\forall~k\in\{0,1,2,\dots\}~t\in(0,\infty).
\end{eqnarray*}
(Note that $\Psi_k$ may take $+\infty$ with a positive probability under $\overline{p}_n(d\xi|x)$. If $\overline{\lambda}(a)\equiv \lambda$, then $\{\sum_{i=0}^n \Psi_i\}_{n\ge 0}$ forms a standard Poisson point process, justifying the use of the prefix ``pseudo'' here.)
\end{definition}

Given a pseudo-Poisson-related policy $\overline{S}^P= \{\overline{p}_n\}_{n=0}^\infty $ and initial state $x_0\in\textbf{X}$, there is a unique probability measure ${\rm P}_{x_0}^{\overline{S}^P}$ on $(\Omega,{\cal F})$ such that ${\rm P}_{x_0}^{\overline{S}^P}(X_0\in dx)=\delta_{x_0}(dx)$,  for each $n\ge 1$ and  $\Gamma_1\in {\cal B}([0,\infty))$, $\Gamma_2\in {\cal B}(\textbf{X})$,
\begin{eqnarray}\label{2021June30Eqn27}
&&{\rm P}_{x_0}^{\overline{S}^P}(  \Theta_n\in \Gamma_1,~X_n\in \Gamma_2|H_{n-1})= \int_{{\bf \Xi}^{GO}}{\rm P}_{x_0}^{\overline{S}^P,\xi}(  \Theta_n\in \Gamma_1,~X_n\in \Gamma_2|H_{n-1}) \overline{p}_n(d\xi|x)\nonumber\\
&:=& \int_{{\bf \Xi}^{GO}}\left\{\int_{\Gamma_1}e^{-\int_0^s q^{GO,\xi}_{X_{n-1}}(t)dt}\tilde{q}^{GO,\xi}(\Gamma_2|{X_{n-1},s)}ds\right\} \overline{p}_n(d\xi|x);\nonumber\\
&&{\rm P}_{x_0}^{\overline{S}^P}(\Theta_n=\infty,~X_n=x_\infty|H_{n-1})=\int_{{\bf \Xi}^{GO}}{\rm P}_{x_0}^{\overline{S}^P,\xi}(  \Theta_n=\infty,~X_n=x_\infty|H_{n-1}) \overline{p}_n(d\xi|x)\nonumber\\
&:=&\int_{{\bf \Xi}^{GO}}\left\{ e^{-\int_0^\infty q^{GO,\xi}_{X_{n-1}}(t)dt}\right\}\overline{p}_n(d\xi|x);
\end{eqnarray}
and
\begin{eqnarray*}
{\rm P}_{x_0}^{\overline{S}^P}(\Theta_n=\infty,~X_n\in \Gamma_2|H_{n-1})= {\rm P}_{x_0}^{\overline{S}^P}(\Theta_n\in\Gamma_1,~X_n=x_\infty|H_{n-1})=0,
\end{eqnarray*}
where
\begin{eqnarray*}
&&q^{GO,\xi}(dy|x,s):=\sum_{k=0}^\infty
q^{GO}(dy|x,\alpha_k)I\{s\in(\tau_k,\tau_{k+1}]\},~\tilde{q}^{GO,\xi}(dy|x,s):=\sum_{k=0}^\infty
\tilde{q}^{GO}(dy|x,\alpha_k)I\{s\in(\tau_k,\tau_{k+1}]\},\\
&&q_x^{GO,\xi}(s):=\sum_{k=0}^\infty
q^{GO}_x(\alpha_k)I\{s\in(\tau_k,\tau_{k+1}]\}.
\end{eqnarray*}
Let the expectation corresponding to ${\rm P}^{\overline{S}^P}_{x_0}$ be denoted as ${\rm E}_{x_0}^{\overline{S}^P}$.


The system performance under $\overline{S}^P$ is measured by
\begin{eqnarray*}
&&{W}_i(x_0,\overline{S}^P):={\rm E}_{x_0}^{\overline{S}^P}\left[\sum_{n=0}^\infty I\{X_n\ne x_\infty\}\int_{{\bf \Xi}^{GO}}\int_{(0,\infty]}  \int_0^t c_i^{GO,\xi}(X_n,s)ds {\rm P}_n^{\overline{S}^P,\xi}(\Theta_{n+1}\in dt|X_n)\overline{p}_n(d\xi|X_n) \right]\\
&:=&{\rm E}_{x_0}^{\overline{S}^P}\left[\sum_{n\ge 0}I\{X_n\ne x_\infty\} \int_{{\bf \Xi}^{GO}}\int_{(0,\infty]}  \int_0^t c_i^{GO+,\xi}(X_n,s)ds {\rm P}_n^{\overline{S}^P,\xi}(\Theta_{n+1}\in dt|X_n)\overline{p}_n(d\xi|X_n)\right]\\
&&-{\rm E}_{x_0}^{\overline{S}^P}\left[\sum_{n\ge 0}I\{X_n\ne x_\infty\} \int_{{\bf \Xi}^{GO}}\int_{(0,\infty]}  \int_0^t c_i^{GO-,\xi}(X_n,s)ds {\rm P}_n^{\overline{S}^P,\xi}(\Theta_{n+1}\in dt|X_n)\overline{p}_n(d\xi|X_n)\right],
\end{eqnarray*}
where $\infty-\infty:=\infty,$ ${\rm P}_n^{\overline{S}^P,\xi}(\Theta_{n+1}\in dt|X_n)$ is defined in (\ref{2021June30Eqn27}), see the terms inside the parentheses therein,
 and \begin{eqnarray*}
c_i^{GO,\xi}(x,s):=\sum_{k=0}^\infty
c_i^{GO}(x,\alpha_k)I\{s\in(\tau_k,\tau_{k+1}]\},~c_i^{GO\pm,\xi}(x,s):=\sum_{k=0}^\infty
c_i^{GO\pm}(x,\alpha_k)I\{s\in(\tau_k,\tau_{k+1}]\}.
\end{eqnarray*}

\begin{theorem}\label{2021June30Theorem05}
Each Markov policy $\overline{S}^M=\{\overline{F}^M_n\}_{n\ge 0}$ in ${\cal M}^{GO}$ can be replicated by a pseudo-Poisson-related policy $\overline{S}^P$ in ${\cal M}^{GO}.$
\end{theorem}

\par\noindent\textit{Proof.} Let some Markov policy $\overline{S}^M=\{\overline{F}^M_n\}_{n\ge 0}$ in ${\cal M}^{GO}$ be given, and define the following $\overline{S}^P=\{\overline{p}_n\}_{n\ge 0}$ by
\begin{eqnarray}\label{2021June30Eqn28}
\overline{p}_{n,0}(da|x):=\int_0^\infty e^{-\int_0^t (\overline{\lambda}+q^{GO}_x)(\overline{F}_n^M,s)ds}(q^{GO}_x(a)+\overline{\lambda}(a))\overline{F}_n^M(x)_t(da)dt
\end{eqnarray}
where
\begin{eqnarray*}
(\overline{\lambda}+q^{GO}_x)(\overline{F}_n^M,s):=\int_{\textbf{A}}(\overline{\lambda}(a)+q^{GO}_x(a))\overline{F}_n^M(x)_s(da);
\end{eqnarray*}
and for each $k\ge 1$,
\begin{eqnarray}\label{2021June30Eqn29}
\overline{p}_{n,k}(da|x)&:=&\int_0^\infty \frac{\overline{\lambda}(\overline{F}_n^M,w) \left(\int_0^w \overline{\lambda}(\overline{F}_n^M,u)du \right)^{k-1}}{(k-1)!}\nonumber\\
     &&  \times\left(\frac{\int_w^\infty e^{-\int_0^t (\overline{\lambda}+q^{GO}_x)(\overline{F}_n^M,s)ds}(q^{GO}_x(a)+\overline{\lambda}(a))\overline{F}_n^M(x)_t(da)dt}{\int_0^\infty \frac{\overline{\lambda}(\overline{F}_n^M,w) \left(\int_0^w \overline{\lambda}(\overline{F}_n^M,u)du \right)^{k-1}}{(k-1)!} e^{-\int_0^w (\overline{\lambda}+q_x^{GO})(\overline{F}_n^M,s)ds }dw}\right)dw
\end{eqnarray}
if the denominator does not vanish, otherwise $\overline{p}_{n,k}(da|x)$ is put to be a fixed probability measure $\overline{p}^\ast(da)$ with $\overline{p}^\ast\in {\cal P}(\textbf{A})$ being concentrated on $\textbf{A}^I$.

For notational convenience, let us introduce
\begin{eqnarray*}
Q_{n,k}(w,x):=\frac{\overline{\lambda}(\overline{F}_n^M,w) \left(\int_0^w \overline{\lambda}(\overline{F}_n^M,u)du \right)^{k-1}}{(k-1)!} e^{-\int_0^w (\overline{\lambda}+q_x^{GO})(\overline{F}_n^M,s)ds },
\end{eqnarray*}
so that
\begin{eqnarray*}
\overline{p}_{n,k}(da|x)&:=&\int_0^\infty \frac{\overline{\lambda}(\overline{F}_n^M,w) \left(\int_0^w \overline{\lambda}(\overline{F}_n^M,u)du \right)^{k-1}}{(k-1)!}\nonumber\\
     &&  \times\left(\frac{\int_w^\infty e^{-\int_0^t (\overline{\lambda}+q^{GO}_x)(\overline{F}_n^M,s)ds}(q^{GO}_x(a)+\overline{\lambda}(a))\overline{F}_n^M(x)_t(da)dt}{\int_0^\infty Q_{n,k}(w,x) dw}\right)dw.
\end{eqnarray*}

It is useful to observe that if
\begin{eqnarray}\label{2021June30Eqn30}
\int_0^\infty Q_{n,k}(w,x)dw:=\int_0^\infty \left\{\frac{\overline{\lambda}(\overline{F}_n^M,w) \left(\int_0^w \overline{\lambda}(\overline{F}_n^M,u)du \right)^{k-1}}{(k-1)!} e^{-\int_0^w (\overline{\lambda}+q_x^{GO})(\overline{F}_n^M,s)ds }\right\}dw
\end{eqnarray}
vanishes
for some $k\ge 1,$ then so does  $\int_0^\infty Q_{n,l}(w,x)dw$ for all $l\in\{1,2,\dots\}.$

First of all, let us verify that
\begin{eqnarray}\label{2021June30Eqn09}
{\rm P}_{x_0}^{\overline{S}^P}(X_n\in dy)={\rm P}_{x_0}^{\overline{S}^M}(X_n \in dy)
\end{eqnarray} as follows.
The case of $n=0$ is evident. Suppose it holds for some $n\ge 0$, and let us prove  ${\rm P}_{x_0}^{\overline{S}^P}(X_{n+1}\in \Gamma|X_n=x)={\rm P}_{x_0}^{\overline{S}^M}(X_{n+1}\in \Gamma|X_n=x)$ for each $x\in \textbf{X}$ and $\Gamma\in{\cal B}(\textbf{X})$, as follows. Note that
\begin{eqnarray}\label{2021June30Eqn06}
&&{\rm P}_{x_0}^{\overline{S}^P}(X_{n+1}\in \Gamma|X_n=x)= \int_{{\bf \Xi}^{GO}}\int_{0}^\infty e^{-\int_0^s q^{GO,\xi}_{x}(t)dt}\tilde{q}^{GO,\xi}(\Gamma|{x,s)}ds \overline{p}_n(d\xi|x)\\
&=&\sum_{k\ge 0}\int_{{\bf \Xi}^{GO}} \int_{(\tau_k,\tau_{k+1})} e^{-\int_0^s q^{GO,\xi}_{x}(t)dt}\tilde{q}^{GO,\xi}(\Gamma|{x,s)}ds \overline{p}_n(d\xi|x)\nonumber\\
&=& \sum_{k=0}^\infty\int_{{\bf \Xi}^{GO}} \int_{(\tau_k,\tau_{k+1})} \tilde{q}^{GO}(\Gamma|x,\alpha_k)\prod_{i=0}^{k-1} e^{-\psi_{i+1} q_x^{GO}(\alpha_i)}e^{-(s-\tau_k)q_x^{GO}(\alpha_k)}ds\overline{p}_n(d\xi|x)\nonumber\\
&=&\sum_{k=0}^\infty\int_{{\bf \Xi}^{GO}}\prod_{i=0}^{k-1} e^{-\psi_{i+1} q_x^{GO}(\alpha_i)} I\{\psi_{i+1}<\infty\} \tilde{q}^{GO}(\Gamma|x,\alpha_k)\int_{0}^{\psi_{k+1}} e^{-q_x^{GO}(\alpha_k)s}ds\overline{p}_n(d\xi|x).\nonumber
\end{eqnarray}
Since $(\Phi_0,\Psi_1),(\Phi_1,\Psi_2),\dots$ are mutually independent under $\overline{p}_n(d\xi|x)$, we see, upon computing the integrals with respect to $\overline{p}_n(d\xi|x)$ in the above, that
\begin{eqnarray}\label{2021June30Eqn02}
&&{\rm P}_{x_0}^{\overline{S}^P}(X_{n+1}\in \Gamma|X_n=x)= \sum_{k=0}^\infty \prod_{i=0}^{k-1}\int_{\textbf{A}}\frac{\overline{\lambda}(a)}{\overline{\lambda}(a)+q^{GO}_x(a)}\overline{p}_{n,i}(da|x)\int_{\textbf{A}} \frac{\tilde{q}^{GO}(\Gamma|x,a)}{\overline{\lambda}(a)+q^{GO}_x(a)}\overline{p}_{n,k}(da|x),
\end{eqnarray}
where we recall that $\overline{\lambda}(a)+q^{GO}_x(a)\ge \min\{1,\lambda\}>0$ for all $a\in\textbf{A}.$
Let us verify for $k\ge 1$ that
\begin{eqnarray}\label{2021June30Eqn01}
\prod_{i=0}^{k-1}\int_{\textbf{A}}\frac{\overline{\lambda}(a)}{\overline{\lambda}(a)+q^{GO}_x(a)}\overline{p}_{n,i}(da|x)=\int_0^\infty Q_{n,k}(w,x)dw
\end{eqnarray}
as follows.
When $k=1,$ the left hand side can be written as
\begin{eqnarray*}
&&\int_{\textbf{A}}\frac{\overline{\lambda}(a)}{\overline{\lambda}(a)+q^{GO}_x(a)}\overline{p}_{n,0}(da|x)\\
&=&\int_{\textbf{A}}\frac{\overline{\lambda}(a)}{\overline{\lambda}(a)+q^{GO}_x(a)}\int_0^\infty e^{-\int_0^t (\overline{\lambda}+q^{GO}_x)(\overline{F}_n^M,s)ds}(q^{GO}_x(a)+\overline{\lambda}(a))\overline{F}_n^M(x)_t(da)dt\\
&=&\int_0^\infty \overline{\lambda}(\overline{F}_n^M,t) e^{-\int_0^t (\overline{\lambda}+q^{GO}_x)(\overline{F}_n^M,s)ds}dt=\int_0^\infty Q_{n,1}(w,x)dw,
\end{eqnarray*}
as desired. (Again, we used here the fact that $\overline{\lambda}(a)+q^{GO}_x(a)\ge \min\{1,\lambda\}>0$ for all $a\in\textbf{A}.$)

Now assume that (\ref{2021June30Eqn01}) holds for some $k\ge 1,$ and we now need show that
\begin{eqnarray*}
\prod_{i=0}^{k}\int_{\textbf{A}}\frac{\overline{\lambda}(a)}{\overline{\lambda}(a)+q^{GO}_x(a)}\overline{p}_{n,i}(da|x)=\int_0^\infty Q_{n,k+1}(w,x)dw.
\end{eqnarray*}
The case when the right hand side vanishes is trivial, because it implies the same for the left hand side  by the definition of $\overline{p}_{n,i}$, see (\ref{2021June30Eqn28}) and (\ref{2021June30Eqn29}), and the observation below (\ref{2021June30Eqn30}).  Thus, we assume that $\int_0^\infty Q_{n,k+1}(w,x)dw>0$, which is equivalent to that $\int_0^\infty Q_{n,k}(w,x)dw>0$ for all $k\ge 1$ as was observed below (\ref{2021June30Eqn30}). Then, by the inductive supposition,
\begin{eqnarray*}
&&\prod_{i=0}^{k}\int_{\textbf{A}}\frac{\overline{\lambda}(a)}{\overline{\lambda}(a)+q^{GO}_x(a)}\overline{p}_{n,i}(da|x)=\prod_{i=0}^{k-1}\int_{\textbf{A}}\frac{\overline{\lambda}(a)}{\overline{\lambda}(a)+q^{GO}_x(a)}\overline{p}_{n,i}(da|x) \int_{\textbf{A}}\frac{\overline{\lambda}(a)}{\overline{\lambda}(a)+q^{GO}_x(a)}\overline{p}_{n,k}(da|x)\\
&=&\int_0^\infty Q_{n,k}(w,x)dw\int_{\textbf{A}}\frac{\overline{\lambda}(a)}{\overline{\lambda}(a)+q^{GO}_x(a)}\overline{p}_{n,k}(da|x).
\end{eqnarray*}
The above expression is equal to
\begin{eqnarray*}
&&\int_0^\infty Q_{n,k}(w,x)dw\int_{\textbf{A}}\frac{\overline{\lambda}(a)}{\overline{\lambda}(a)+q^{GO}_x(a)}\int_0^\infty \frac{\overline{\lambda}(\overline{F}_n^M,w) \left(\int_0^w \overline{\lambda}(\overline{F}_n^M,u)du \right)^{k-1}}{(k-1)!}\\
     &&  \times\left(\frac{\int_w^\infty e^{-\int_0^t (\overline{\lambda}+q^{GO}_x)(\overline{F}_n^M,s)ds}(q^{GO}_x(a)+\overline{\lambda}(a))\overline{F}_n^M(x)_t(da)dt}{\int_0^\infty Q_{n,k}(w,x)dw}\right)dw\\
&=&\int_{\textbf{A}}\frac{\overline{\lambda}(a)}{\overline{\lambda}(a)+q^{GO}_x(a)}\int_0^\infty \frac{\overline{\lambda}(\overline{F}_n^M,w) \left(\int_0^w \overline{\lambda}(\overline{F}_n^M,u)du \right)^{k-1}}{(k-1)!}\\
     &&  \times\left( \int_w^\infty e^{-\int_0^t (\overline{\lambda}+q^{GO}_x)(\overline{F}_n^M,s)ds}(q^{GO}_x(a)+\overline{\lambda}(a))\overline{F}_n^M(x)_t(da)dt \right)dw\\
&=&\int_0^\infty \left[\int_w^\infty   \overline{\lambda}(\overline{F}_n^M,t) e^{-\int_0^t (\overline{\lambda}+q^{GO}_x)(\overline{F}_n^M,s)ds}  dt\right] \frac{\overline{\lambda}(\overline{F}_n^M,w) \left(\int_0^w \overline{\lambda}(\overline{F}_n^M,u)du \right)^{k-1}}{(k-1)!} dw.
\end{eqnarray*}
Integrating by parts the above integral, we may write the previous expression as
\begin{eqnarray*}
&&\left[\int_w^\infty  \overline{\lambda}(\overline{F}_n^M,t) e^{-\int_0^t (\overline{\lambda}+q^{GO}_x)(\overline{F}_n^M,s)ds}  dt \frac{  \left(\int_0^w \overline{\lambda}(\overline{F}_n^M,u)du \right)^k}{k!} \right]_{0}^\infty\\
&&+\int_0^\infty  \frac{  \left(\int_0^w \overline{\lambda}(\overline{F}_n^M,u)du \right)^k}{k!} \overline{\lambda}(\overline{F}_n^M,w) e^{-\int_0^w (\overline{\lambda}+q^{GO}_x)(\overline{F}_n^M,s)ds}dw=\int_0^\infty Q_{n,k+1}(w,x)dw,
\end{eqnarray*}
where for the equality, one may apply routine analysis based on $\overline{\lambda}(a)+q^{GO}_x(a)\ge \min\{1,\lambda\}>0$ for all $a\in\textbf{A}.$ This thus proves (\ref{2021June30Eqn01}) for all $k\ge 1$.

We may substitute (\ref{2021June30Eqn01}) back in (\ref{2021June30Eqn02}):
\begin{eqnarray}\label{2021June30Eqn05}
&&{\rm P}_{x_0}^{\overline{S}^P}(X_{n+1}\in \Gamma|X_n=x)= \int_{\textbf{A}} \frac{\tilde{q}^{GO}(\Gamma|x,a)}{\overline{\lambda}(a)+q^{GO}_x(a)}\overline{p}_{n,0}(da|x)\nonumber\\
&&+\sum_{k=1}^\infty \int_0^\infty Q_{n,k}(w,x)dw\int_{\textbf{A}} \frac{\tilde{q}^{GO}(\Gamma|x,a)}{\overline{\lambda}(a)+q^{GO}_x(a)}\overline{p}_{n,k}(da|x)\nonumber\\
&=&  \int_0^\infty \tilde{q}^{GO}(\Gamma|x,\overline{F}_n^M,t)  e^{-\int_0^t (\overline{\lambda}+q^{GO}_x)(\overline{F}_n^M,s)ds}dt\nonumber\\
&&+ \sum_{k=1}^\infty   \int_0^\infty \frac{\overline{\lambda}(\overline{F}_n^M,w) \left(\int_0^w \overline{\lambda}(\overline{F}_n^M,u)du \right)^{k-1}}{(k-1)!} \tilde{q}^{GO}(\Gamma|x,\overline{F}_n^M,t) \int_w^\infty e^{-\int_0^t (\overline{\lambda}+q^{GO}_x)(\overline{F}_n^M,s)ds}dtdw,
\end{eqnarray}
with the above equalities being valid no matter $\int_0^\infty Q_{n,k}(w,x)dw$ vanishes or not: indeed, if
\begin{eqnarray*}
\int_0^\infty Q_{n,k}(w,x)dw=0,
\end{eqnarray*}
then the summands in the last one of the previous equalities vanish, too.

Note that
\begin{eqnarray*}
&&\sum_{k=1}^\infty   \int_0^\infty \frac{\overline{\lambda}(\overline{F}_n^M,w) \left(\int_0^w \overline{\lambda}(\overline{F}_n^M,u)du \right)^{k-1}}{(k-1)!} \tilde{q}^{GO}(\Gamma|x,\overline{F}_n^M,t) \int_w^\infty e^{-\int_0^t (\overline{\lambda}+q^{GO}_x)(\overline{F}_n^M,s)ds}dtdw\\
&=&\sum_{k=1}^\infty   \int_0^\infty \int_0^t  \overline{\lambda}(\overline{F}_n^M,w) \left\{\frac{ \left(\int_0^w \overline{\lambda}(\overline{F}_n^M,u)du \right)^{k-1}}{(k-1)!}\right\} dw\tilde{q}^{GO}(\Gamma|x,\overline{F}_n^M,t)   e^{-\int_0^t (\overline{\lambda}+q^{GO}_x)(\overline{F}_n^M,s)ds}dt \\
&=&\sum_{k=1}^\infty   \int_0^\infty \int_0^t  \overline{\lambda}(\overline{F}_n^M,w) \int_{\{0\le v_1\le v_2\le\dots\le v_{k-1}\le w\}} \prod_{j=1}^{k-1}\overline{\lambda}(\overline{F}_n^M,v_j)dv_1dv_2\dots dv_{k-1} dw\\
&&\times \tilde{q}^{GO}(\Gamma|x,\overline{F}_n^M,t)   e^{-\int_0^t (\overline{\lambda}+q^{GO}_x)(\overline{F}_n^M,s)ds}dt\\
&=&\sum_{k=1}^\infty   \int_0^\infty    \int_{\{0\le v_1\le v_2\le\dots\le v_{k-1}\le w\le t\}} \prod_{j=1}^{k-1}\overline{\lambda}(\overline{F}_n^M,v_j)  \overline{\lambda}(\overline{F}_n^M,w) dv_1dv_2\dots dv_{k-1} dw\\
&&\times \tilde{q}^{GO}(\Gamma|x,\overline{F}_n^M,t)   e^{-\int_0^t (\overline{\lambda}+q^{GO}_x)(\overline{F}_n^M,s)ds}dt\\
&=&\sum_{k=1}^\infty   \int_0^\infty    \left\{\frac{ \left(\int_0^t \overline{\lambda}(\overline{F}_n^M,u)du \right)^{k}}{k!}\right\}\tilde{q}^{GO}(\Gamma|x,\overline{F}_n^M,t)   e^{-\int_0^t (\overline{\lambda}+q^{GO}_x)(\overline{F}_n^M,s)ds}dt
\end{eqnarray*}
where the first equality is by the Fubini-Tonelli theorem, and for the second as well as the last equality, recall the following equality, which is valid for any real-valued integrable function $f$:
\begin{eqnarray*}
\left(\int_0^w f(u)du \right)^{k-1}=(k-1)!\int_{\{0\le v_1\le v_2\le\dots\le v_{k-1}\le w\}}\prod_{j=1}^{k-1} f(v_j)dv_1dv_2\dots d v_{k-1}.
\end{eqnarray*}
 With the above equalities, (\ref{2021June30Eqn05}) can be written as follows:
\begin{eqnarray}\label{2021June30Eqn07}
&&{\rm P}_{x_0}^{\overline{S}^P}(X_{n+1}\in \Gamma|X_n=x)=\sum_{k=0}^\infty   \int_0^\infty    \left\{\frac{ \left(\int_0^t \overline{\lambda}(\overline{F}_n^M,u)du \right)^{k}}{k!}\right\}\tilde{q}^{GO}(\Gamma|x,\overline{F}_n^M,t)   e^{-\int_0^t (\overline{\lambda}+q^{GO}_x)(\overline{F}_n^M,s)ds}dt\nonumber\\
&=& \int_0^\infty     \tilde{q}^{GO}(\Gamma|x,\overline{F}_n^M,t)   e^{-\int_0^t q^{GO}_x(\overline{F}_n^M,s)ds}dt\\
&=&{\rm P}_{x_0}^{\overline{S}^M}(X_{n+1}\in \Gamma|X_n=x),\nonumber
\end{eqnarray}
as desired.

The rest verifies
\begin{eqnarray}\label{2021June30Eqn08}
&&{\rm E}_{x_0}^{\overline{S}^P}\left[I\{X_n\ne x_\infty\} \int_{{\bf \Xi}^{GO}}\int_{(0,\infty]}  \int_0^t c_i^{GO,\xi}(X_n,s)ds {\rm P}_n^{\overline{S}^P,\xi}(\Theta_{n+1}\in dt|X_n)\overline{p}_n(d\xi|X_n)\right]\nonumber\\
&=&{\rm E}_{x_0}^{\overline{S}^M}\left[I\{X_n\ne x_\infty\}\int_0^{\Theta_{n+1}} c_i^{GO}(X_n,\overline{F}_n^M,s)ds\right],
\end{eqnarray}
which would complete the proof of this theorem.
It is sufficient to assume in the rest of this proof that $c_i^{GO}$ is nonnegative and bounded on $\textbf{X}\times \textbf{A}$: the general case can be handled based on this simpler case with the help of the monotone convergence theorem.

Note that on $\{X_n\ne x_\infty\},$
\begin{eqnarray*}
\int_{(0,\infty]}  \int_0^t c_i^{GO,\xi}(X_n,s)ds {\rm P}_n^{\overline{S}^P,\xi}(\Theta_{n+1}\in dt|X_n)={\rm E}_n^{\overline{S}^P,\xi}\left[ \int_0^{\Theta_{n+1}} c_i^{GO,\xi}(X_n,s)ds |X_n\right],
\end{eqnarray*}
where ${\rm E}_n^{\overline{S}^P,\xi}[\cdot|X_n]$ is understood with respect to ${\rm P}_n^{\overline{S}^P,\xi}(\Theta_{n+1}\in dt|X_n),$ which is defined in (\ref{2021June30Eqn27}), see the terms inside the parentheses therein.

Now, the left hand side of (\ref{2021June30Eqn08}) can be written as
\begin{eqnarray*}
&&{\rm E}_{x_0}^{\overline{S}^P}\left[I\{X_n\ne x_\infty\} \int_{{\bf \Xi}^{GO}}\int_{(0,\infty]}  \int_0^t c_i^{GO,\xi}(X_n,s)ds {\rm P}_n^{\overline{S}^P,\xi}(\Theta_{n+1}\in dt|X_n)\overline{p}_n(d\xi|X_n)\right]\nonumber\\
&=&{\rm E}_{x_0}^{\overline{S}^P}\left[I\{X_n\ne x_\infty\}\int_{{\bf\Xi}^{GO}}{\rm E}^{\overline{S}^P,\xi}_{x_0}\left[\int_0^{\infty} c_i^{GO,\xi}(X_n,s) I\{s<\Theta_{n+1}\}ds|X_n\right] \overline{p}_n(d\xi|X_n)\right]\nonumber\\
&=&\int_{\textbf{X}}  {\rm P}^{\overline{S}^P}_{x_0}(X_n\in dx) \left\{ \int_{{\bf\Xi}^{GO}} \int_0^\infty  c_i^{GO,\xi}(x,s)e^{-\int_0^s q^{GO,\xi}_{x}(t)dt} ds \overline{p}_n(d\xi|x)\right\}.
\end{eqnarray*}
Note that the term inside the parenthesis is in the same form as the term on the right hand side of (\ref{2021June30Eqn06}), where $\tilde{q}^{GO,\xi}(\Gamma|{x,s)}$ is replaced by  $c_i^{GO,\xi}(x,s)$ with the latter term having been assumed to be nonnegative and bounded. Therefore, the calculations in (\ref{2021June30Eqn06})-(\ref{2021June30Eqn07}) apply with obvious modifications (more precisely, replacing $\tilde{q}^{GO}(\Gamma|x,a)$ by $c_i^{GO}(x,a)$), leading to
\begin{eqnarray}\label{2021June30Eqn17}
&&{\rm E}_{x_0}^{\overline{S}^P}\left[I\{X_n\ne x_\infty\} \int_{{\bf \Xi}^{GO}}\int_{(0,\infty]}  \int_0^t c_i^{GO,\xi}(X_n,s)ds {\rm P}_n^{\overline{S}^P,\xi}(\Theta_{n+1}\in dt|X_n)\overline{p}_n(d\xi|X_n)\right]\nonumber\\
&=&\int_{\textbf{X}}{\rm P}^{\overline{S}^P}_{x_0}(X_n\in dx)   \left\{ \sum_{k=0}^\infty \prod_{i=0}^{k-1}\int_{\textbf{A}}\frac{\overline{\lambda}(a)}{\overline{\lambda}(a)+q^{GO}_x(a)}\overline{p}_{n,i}(da|x)\int_{\textbf{A}} \frac{ c_i^{GO}(x,a)}{\overline{\lambda}(a)+q^{GO}_x(a)}\overline{p}_{n,k}(da|x)\right\}\\
&=&
\int_{\textbf{X}}{\rm P}^{\overline{S}^P}_{x_0}(X_n\in dx)   \left\{ \int_0^\infty     c_i^{GO}(x,\overline{F}_n^M,t)   e^{-\int_0^t q^{GO}_x(\overline{F}_n^M,s)ds}dt\right\},\nonumber
\end{eqnarray}
where for the first and the second equality, compare the corresponding terms in the parentheses with (\ref{2021June30Eqn02}) and (\ref{2021June30Eqn07}).

On the other hand, the right hand side of (\ref{2021June30Eqn08}) can be written as
\begin{eqnarray*}
&&{\rm E}_{x_0}^{\overline{S}^M}\left[I\{X_n\ne x_\infty\}{\rm E}_{x_0}^{\overline{S}^M}\left[ \int_0^{\infty} c_i^{GO}(X_n,\overline{F}_n^M,s)I\{s<\Theta_{n+1}\}ds|X_n\right]\right]\\
&=& \int_{\textbf{X}}{\rm P}^{\overline{S}^M}_{x_0}(X_n\in dx)   \left\{ \int_0^\infty     c_i^{GO}(x,\overline{F}_n^M,t)   e^{-\int_0^t q^{GO}_x(\overline{F}_n^M,s)ds}dt\right\}.
\end{eqnarray*}
Since ${\rm P}^{\overline{S}^M}_{x_0}(X_n\in dx)={\rm P}^{\overline{S}^P}_{x_0}(X_n\in dx)$ as was verified earlier in this proof, c.f., (\ref{2021June30Eqn09}), we see that the previous expression coincides with the term on the left hand side of (\ref{2021June30Eqn08}), as required. $\hfill\Box$
\bigskip

\subsection{Poisson-related strategy in the gradual-impulsive control model ${\cal M}$}

Recall that $\lambda\in (0,\infty)$ is a fixed constant.
Let ${\bf \Xi}:=[0,\infty)\times \textbf{A}^G\times ((0,\infty)\times \textbf{A}^G)^\infty$ be the countable product.
The context should exclude any confusion that the generic notation for an element of ${\bf \Xi}^{GO}$ is still $\xi=\{(\psi_n,\alpha_n)\}_{n\ge 0}\in {\bf \Xi}$, and the coordinate random variables are still denoted, for each $\xi=\{(\psi_n,\alpha_n)\}_{n\ge 0}\in {\bf \Xi}$, by $\Psi_n(\xi):=\psi_n$ and $\Phi_n(\xi):=\alpha_n$.
For each $n\in \{0,1,\dots\},$ let $p_n(d\xi|x)$ be a stochastic kernel on ${\cal B}({\bf \Xi})$ given $x\in\textbf{X}$, which is specified by the following: for each $x\in \textbf{X},$  under $p_n(d\xi|x)$, the coordinate random variables $ \Psi_0,\Phi_0,\Psi_1,\Phi_1,\dots$ are mutually independent and
\begin{eqnarray*}
&&p_n(\Psi_0\in dt|x)=\delta_0(dt),~p_n(\Psi_k\le t|x)=1-e^{-\lambda t},~\forall~k\in\{1,2,\dots\},\\
&&p_n(\Phi_k\in da|x)=:p_{n,k}(da|x)~\forall~k\in\{0, 1,2,\dots\}.
\end{eqnarray*}
(Hence, under $p_n(d\xi|x)$, $\{\sum_{k=0}^n \Psi_k\}_{n\ge 1}$ is a Poisson point process.)
Let $\sigma_n^{P,(0)}(d\hat{c}\times d\hat{b}|x,\xi)$ be a stochastic kernel on ${\cal B}([0,\infty]\times \textbf{A}^I)$ from $(x,\xi)\in\textbf{X}\times {\bf \Xi}$.

\begin{definition}[Poisson-related strategy for ${\cal M}$]
The pairs $\{(\sigma_n^{P,(0)},p_n)\}_{n\ge 0}=:\sigma^P$ is called a Poisson-related strategy in the gradual-impulsive control model ${\cal M}$.
\end{definition}

Given $\xi=\{(\psi_n,\alpha_n)\}_{n\ge 0}\in {\bf \Xi}$, with the generic notation $\tau_n:=\sum_{k=0}^n \psi_k$ for each $n\in\{0,1,\dots\},$ we put
\begin{eqnarray}\label{2021June30Eqn32}
&&q^\xi(dy|x,s):=\sum_{k=0}^\infty
q(dy|x,\alpha_k)I\{s\in(\tau_k,\tau_{k+1}]\},~\tilde{q}^\xi(dy|x,s):=\sum_{k=0}^\infty
\tilde{q}(dy|x,\alpha_k)I\{s\in(\tau_k,\tau_{k+1}]\},\nonumber\\
&&q_x^\xi(s):=\sum_{k=0}^\infty
q_x(\alpha_k)I\{s\in(\tau_k,\tau_{k+1}]\}.
\end{eqnarray}

Under a Poisson-related strategy $\sigma^P=\{(\sigma_n^{P,(0)},p_n)\}_{n\ge 0}$, the transition law of $\hat{X}_{n+1}=(\hat{\Theta}_{n+1},X_{n+1} )$ given $(\hat{\Theta}_n,X_n)=(\theta,x)\in \hat{\textbf{X}}$ is denoted by $G_n^{\sigma^P}$, which is defined for each bounded measurable function $g$ on $\hat{\textbf{X}}$ by
\begin{eqnarray}\label{2021June30Eqn11}
&&\int_{\hat{\textbf{X}}}g(t,y)G_n^{\sigma^P}(dt\times dy|(\theta,x))\nonumber\\
&:=&\int_{[0,\infty]\times \textbf{A}^I \times{\bf \Xi}}\left\{ \int_0^{\hat{c}} \int_{\textbf{X}} g(t,y)\tilde{q}^\xi(dy|x,t)e^{-\int_0^t q^\xi_x(s)ds}dt+I\{\hat{c}=\infty\} g(\infty,x_\infty)e^{-\int_0^\infty q^\xi_x(s)ds}\right.\nonumber\\
&&\left.+I\{\hat{c}<\infty\}e^{-\int_0^{\hat{c}} q^\xi_x(s)ds}\int_{\textbf{X}}g(\hat{c},y)Q(dy|x,\hat{b})\right\}\sigma_n^{P,(0)}(d\hat{c}\times d\hat{b}|x,\xi)p_n(d\xi|x)\nonumber\\
&=:&\int_{[0,\infty]\times \textbf{A}^I \times{\bf \Xi}} \left\{\int_{\hat{\textbf{X}}}g(t,y) G_n^{\sigma^P,\xi}(dt\times dy|(\theta,x),\hat{c},\hat{b})\right\} \sigma_n^{P,(0)}(d\hat{c}\times d\hat{b}|x,\xi)p_n(d\xi|x)
\end{eqnarray}
for each
  $(\theta,x)\in[0,\infty)\times \textbf{X}$; and
$\int_{\hat{\textbf{X}}}g(t,y)G_n^{\sigma^P}(dt\times dy|(\infty,x_\infty)):=g(\infty,x_\infty).
$

\begin{remark}\label{2021June30Remark01}
Note that, $G^{\sigma^P}(dt\times dy|(\theta,x))$ and $G^{\sigma^P,\xi}(dt\times dy|(\theta,x),\hat{c},\hat{b}),$ which is defined in (\ref{2021June30Eqn11}), see the terms inside the parentheses therein, depend on $(\theta,x)$ only through $x\in\textbf{X}_\infty:=\textbf{X}\cup\{x_\infty\}$, and therefore, we will write $G^{\sigma^P}(dt\times dy|x)$ and $G^{\sigma^P,\xi}(dt\times dy|x,\hat{c},\hat{b})$ for $G^{\sigma^P}(dt\times dy|(\theta,x))$ and $G^{\sigma^P,\xi}(dt\times dy|(\theta,x),\hat{c},\hat{b})$ in what follows. The same applies to $l_i^{\sigma^P,n}(\hat{x})=l_i^{\sigma^P,n}(x)$ introduced below.
\end{remark}

The sequence $\{G_n^{\sigma^P}\}_{n\ge 0}$ together with the initial distribution $\delta_{x_0}(dy)\delta_0(dt)$ defines a probability $\hat{\rm P}_{x_0}^{\sigma^P}$ on $\left[\bigcup_{n\ge 1}([0,\infty)\times \textbf{X})^n\times \{(\infty,x_\infty)\}^\infty\right] \cup([0,\infty)\times \textbf{X})^\infty$.
Let ${\hat{\rm E}_{x_0}^{\sigma^P}}$ be the expectation with respect to ${\hat{\rm P}_{x_0}^{\sigma^P}}.$
The system performance under $\sigma^P$ is measured by
\begin{eqnarray*}
\hat{W}_i(x_0,\sigma^P):=\sum_{n\ge 0}\hat{\rm E}_{x_0}^{\sigma^P}\left[  l_i^{\sigma^P,n}(X_n)\right]:= \hat{\rm E}_{x_0}^{\sigma^P}\left[ \sum_{n\ge 0} l_i^{\sigma^P,n(+)}(X_n)\right]- \hat{\rm E}_{x_0}^{\sigma^P}\left[\sum_{n\ge 0} l_i^{\sigma^P,n(-)}(X_n)\right],
\end{eqnarray*}
where we recall the generic notation $\hat{X}_n=(\hat{\Theta}_n,X_n)$ for a state variable in the gradual-impulsive control model ${\cal M},$
and
\begin{eqnarray*}
&&l_i^{\sigma^P,n(\pm)}(x):=\int_{[0,\infty]\times \textbf{A}^I \times{\bf \Xi}}  \int_{0}^\infty I\{x\in \textbf{X}\} \left\{ \int_0^t c_i^{G\pm,\xi}(x,s)ds+I\{t=\hat{c}<\infty\}c_i^{I\pm}(x,\hat{b})\right\}\\
&&G_n^{\sigma^P,\xi}(dt\times \textbf{X}_\infty|x,\hat{c},\hat{b})\sigma_n^{P,(0)}(d\hat{c}\times d\hat{b}|x,\xi)p_n(d\xi|x)~\forall~\hat{x}=(\theta,x)\in\hat{\textbf{X}}
\end{eqnarray*}
with
$c_i^{G\pm,\xi}(x,s):=\sum_{k=0}^\infty
c_i^{G\pm}(x,\alpha_k)I\{s\in(\tau_k,\tau_{k+1}]\}$, and $c_i^{I,\pm}$ being the positive part and the negative part of $c_i^I$, respectively.

If $c_i^I$ and $c_i^G$ are $[0,\infty]$-valued, then
the cost function under $\sigma^P$ over the corresponding sojourn time is given by
\begin{eqnarray}\label{2021June30Eqn37}
&&l_i^{\sigma^P,n}(x):=\int_{[0,\infty]\times \textbf{A}^I \times{\bf \Xi}} \int_{\hat{\textbf{X}}} I\{x\in \textbf{X}\} \left\{ \int_0^t c_i^{G,\xi}(x,s)ds+I\{t=\hat{c}<\infty\}c_i^I(x,\hat{b})\right\}\nonumber\\
&&\times G_n^{\sigma^P,\xi}(dt\times dy|x,\hat{c},\hat{b})\sigma_n^{P,(0)}(d\hat{c}\times d\hat{b}|x,\xi)p_n(d\xi|x)\nonumber\\
&=&\int_{[0,\infty]\times \textbf{A}^I \times{\bf \Xi}}  \int_{0}^\infty I\{x\in \textbf{X}\} \left\{ \int_0^t c_i^{G,\xi}(x,s)ds+I\{t=\hat{c}<\infty\}c_i^I(x,\hat{b})\right\}\nonumber\\
&&\times G_n^{\sigma^P,\xi}(dt\times \textbf{X}_\infty|x,\hat{c},\hat{b})\sigma_n^{P,(0)}(d\hat{c}\times d\hat{b}|x,\xi)p_n(d\xi|x)~\forall~\hat{x}=(\theta,x)\in\hat{\textbf{X}}
\end{eqnarray}
where $\textbf{X}_\infty=\textbf{X}\cup\{x_\infty\}$, and
\begin{eqnarray*}
c_i^{G,\xi}(x,s):=\sum_{k=0}^\infty
c_i^G(x,\alpha_k)I\{s\in(\tau_k,\tau_{k+1}]\}.
\end{eqnarray*}

\begin{theorem}\label{2021June30Theorem02}
Each pseudo-Poisson-related policy $\overline{S}^P=\{\overline{p}_n\}_{n\ge 0}$ in the gradual control model ${\cal M}^{GO}$ can be replicated by a Poisson-related strategy $\sigma^P=\{(\sigma_n^{P,(0)},p_n)\}_{n\ge 0}$ in the gradual-impulsive control model ${\cal M}.$
\end{theorem}

\par\noindent\textit{Proof.} Let a pseudo-Poisson-related policy $\overline{S}^P=\{\overline{p}_n\}_{n\ge 0}$ in the gradual control model ${\cal M}^{GO}$ be fixed. Consider the Poisson-related strategy $\sigma^P=\{(\sigma_n^{P,(0)},p_n)\}_{n\ge 0}$ in the gradual-impulsive control model ${\cal M}$ defined by the following:
on ${\cal B}(\textbf{A}^G),$ for each $x\in \textbf{X}$,
\begin{eqnarray}\label{2021June30Eqn35}
p_{n,k}(da|x):=
\begin{cases}
\frac{\overline{p}_{n,k}(da|x)}{\overline{p}_{n,k}(\textbf{A}^G|x)} & \mbox{if~} \overline{p}_{n,k}(\textbf{A}^G|x)>0;\\
p^\ast(da)&\mbox{otherwise,}
\end{cases}
\end{eqnarray}
where $p^\ast\in {\cal P}(\textbf{A}^G)$ is a fixed probability measure; for each  $x\in\textbf{X}$ and $\xi=(\psi_0,\alpha_0,\psi_1,\alpha_1,\dots)\in {\bf \Xi}$ with $\tau_n=\sum_{k=0}^n \psi_k,$
\begin{eqnarray}\label{2021June30Eqn36}
\sigma_n^{P,(0)}(d\hat{c}\times d\hat{b}|x,\xi):= \sum_{k=0}^\infty \delta_{\tau_k}(d\hat{c}) \overline{p}_{n,k}(d\hat{b}|x)\prod_{m=0}^{k-1} \overline{p}_{n,m}(\textbf{A}^G|x)+\delta_\infty(d\hat{c})\prod_{m=0}^\infty \overline{p}_{n,m}(\textbf{A}^G|x)p^{**}(d\hat{b}),
\end{eqnarray}
where $p^{**}\in{\cal P}(\textbf{A}^I)$ is a fixed probability measure.
Observe that $\sigma_n^{P,(0)}(d\hat{c}\times d\hat{b}|x,\xi)$ defined above depends on $\xi\in{\bf \Xi}$ only through $\xi^-:=(\psi_0,\psi_1,\psi_2,\dots).$

In what follows, we will show in two steps that $\sigma^P$ defined above is a required replicating strategy.

\textit{Step 1.} Firstly, let us verify that
\begin{eqnarray}\label{2021June30Eqn16}
\hat{\rm P}_{x_0}^{\sigma^P}(X_n \in dy)= {\rm P}_{x_0}^{\overline{S}^P}(X_n\in dy).
\end{eqnarray}
Since the above is clearly valid when $n=0$, both sides being equal to $\delta_{x_0}(dy),$ using an inductive argument, it is sufficient to
verify that for an arbitrarily fixed $\Gamma\in{\cal B}(\textbf{X})$ and $x\in \textbf{X}$, for all $n\ge 0$,
\begin{eqnarray}\label{2021June30Eqn10}
G^{\sigma^P}_n([0,\infty)\times \Gamma|x)={\rm P}_{x_0}^{\overline{S}^P}(X_{n+1}\in \Gamma|X_n=x),
\end{eqnarray}
as follows. (Recall (\ref{2021June30Eqn11}) and Remark \ref{2021June30Remark01} for the definition of $G^{\sigma^P,\xi}(dt\times dy|(x,\hat{c},\hat{b})$ with a generic $\sigma^P$.)

Recall that the right hand side of (\ref{2021June30Eqn10}) was computed in (\ref{2021June30Eqn02}), which can be now written out more explicitly using $\textbf{A}=\textbf{A}^I\cup\textbf{A}^G$, $\textbf{A}^I\cap \textbf{A}^G=\emptyset$, $q^{GO}_x(a)=q_x(a)I\{a\in\textbf{A}^G\}+I\{a\in \textbf{A}^I\}$, $\tilde{q}^{GO}(\Gamma|x,a)=Q(\Gamma|x,a)$ for each $a\in\textbf{A}^I$, and $\overline{\lambda}(a)=\lambda I\{a\in\textbf{A}^G\}$ on $\textbf{A}$:
\begin{eqnarray}\label{2021June30Eqn12}
&&{\rm P}_{x_0}^{\overline{S}^P}(X_{n+1}\in \Gamma|X_n=x)= \sum_{k=0}^\infty \left[\prod_{i=0}^{k-1}\int_{\textbf{A}}\frac{\overline{\lambda}(a)}{\overline{\lambda}(a)+q^{GO}_x(a)}\overline{p}_{n,i}(da|x)\right]\int_{\textbf{A}} \frac{\tilde{q}^{GO}(\Gamma|x,a)}{\overline{\lambda}(a)+q^{GO}_x(a)}\overline{p}_{n,k}(da|x)\nonumber\\
&=&\sum_{k=0}^\infty \left[\prod_{i=0}^{k-1}\int_{\textbf{A}^G}\frac{\lambda}{\lambda+q_x(a)}\overline{p}_{n,i}(da|x)\right]\left(\int_{\textbf{A}^G} \frac{\tilde{q}(\Gamma|x,a)}{\lambda+q_x(a)}\overline{p}_{n,k}(da|x)+\int_{\textbf{A}^I} Q(\Gamma|x,a)\overline{p}_{n,k}(da|x)\right)\nonumber\\
&=&\sum_{k=0}^\infty \left[\prod_{i=0}^{k-1}\int_{\textbf{A}^G}\frac{\lambda}{\lambda+q_x(a)}\overline{p}_{n,i}(da|x)\right]\int_{\textbf{A}^G} \frac{\tilde{q}(\Gamma|x,a)}{\lambda+q_x(a)}\overline{p}_{n,k}(da|x) \nonumber\\
&&+\sum_{k=0}^\infty \left[\prod_{i=0}^{k-1}\int_{\textbf{A}^G}\frac{\lambda}{\lambda+q_x(a)}\overline{p}_{n,i}(da|x)\right] \int_{\textbf{A}^I} Q(\Gamma|x,a)\overline{p}_{n,k}(da|x)\nonumber\\
&=:&B_1+B_2.
\end{eqnarray}

On the other hand, the left hand side of (\ref{2021June30Eqn10}) may be written as
\begin{eqnarray}\label{2021June30Eqn19}
&&G^{\sigma^P}_n([0,\infty)\times \Gamma|x)=\int_{{\bf \Xi}} p_n(d\xi|x) \int_{[0,\infty]\times \textbf{A}^I}  \left\{ \int_0^{\hat{c}}  \tilde{q}^\xi(\Gamma|x,t)e^{-\int_0^t q^\xi_x(s)ds}dt\right.\nonumber\\
&&\left.+I\{\hat{c}<\infty\}e^{-\int_0^{\hat{c}} q^\xi_x(s)ds}Q(\Gamma|x,\hat{b})\right\}\sigma_n^{P,(0)}(d\hat{c}\times d\hat{b}|x,\xi)\nonumber\\
&=& \int_{{\bf \Xi}} p_n(d\xi|x) \left(\sum_{k=0}^\infty \int_0^{\tau_k} \tilde{q}^\xi(\Gamma|x,t) e^{-\int_0^t q^\xi_x(s)ds}dt \overline{p}_{n,k}(\textbf{A}^I|x)\prod_{m=0}^{k-1}\overline{p}_{n,m}(\textbf{A}^G|x)\right.\nonumber\\
&&\left.+\int_0^\infty \tilde{q}^\xi(\Gamma|x,t)e^{-\int_0^t q^\xi_x(s)ds}dt  \prod_{m=0}^{\infty}\overline{p}_{n,m}(\textbf{A}^G|x) \right.\nonumber\\
&&\left.+ \sum_{k=0}^\infty e^{-\int_0^{\tau_k} q^\xi_x(s)ds} \int_{\textbf{A}^I}Q(\Gamma|x,\hat{b}) \overline{p}_{n,k}(d\hat{b}|x)\prod_{m=0}^{k-1}\overline{p}_{n,m}(\textbf{A}^G|x)\right),
\end{eqnarray}
where the first equality is by (\ref{2021June30Eqn11}), and the second equality is by the above definition of $\sigma^{P,(0)}_n$, see (\ref{2021June30Eqn36}).
Thus,
\begin{eqnarray}\label{2021June30Eqn15}
&&G^{\sigma^P}_n([0,\infty)\times \Gamma|x)\\
&=& \int_{{\bf \Xi}} p_n(d\xi|x)  \sum_{k=0}^\infty \int_0^{\tau_k} \tilde{q}^\xi(\Gamma|x,t) e^{-\int_0^t q^\xi_x(s)ds}dt \overline{p}_{n,k}(\textbf{A}^I|x)\prod_{m=0}^{k-1}\overline{p}_{n,m}(\textbf{A}^G|x) \nonumber\\
&& +\int_{{\bf \Xi}} p_n(d\xi|x) \int_0^\infty \tilde{q}^\xi(\Gamma|x,t)e^{-\int_0^t q^\xi_x(s)ds}dt  \prod_{m=0}^{\infty}\overline{p}_{n,m}(\textbf{A}^G|x) \nonumber\\
&& + \int_{{\bf \Xi}} p_n(d\xi|x)\sum_{k=0}^\infty e^{-\int_0^{\tau_k} q^\xi_x(s)ds} \int_{\textbf{A}^I}Q(\Gamma|x,\hat{b}) \overline{p}_{n,k}(d\hat{b}|x)\prod_{m=0}^{k-1}\overline{p}_{n,m}(\textbf{A}^G|x)=:C_1+C_2+C_3. \nonumber
\end{eqnarray}
We analyze the above summands term by term as follows.

As for $C_3,$ we see
\begin{eqnarray*}
&&C_3:=\int_{{\bf \Xi}} p_n(d\xi|x)\sum_{k=0}^\infty e^{-\int_0^{\tau_k} q^\xi_x(s)ds} \int_{\textbf{A}^I}Q(\Gamma|x,\hat{b}) \overline{p}_{n,k}(d\hat{b}|x)\prod_{m=0}^{k-1}\overline{p}_{n,m}(\textbf{A}^G|x)\\
&=&\sum_{k=0}^\infty \int_{{\bf \Xi}} p_n(d\xi|x)\left(\left[\prod_{l=0}^{k-1} e^{-\psi_{l+1}q_x(\alpha_{l})}\right] \int_{\textbf{A}^I}Q(\Gamma|x,\hat{b}) \overline{p}_{n,k}(d\hat{b}|x)\prod_{m=0}^{k-1}\overline{p}_{n,m}(\textbf{A}^G|x)\right)\\
&=&\sum_{k=0}^\infty  \left[\prod_{l=0}^{k-1} \int_{\textbf{A}^G} \frac{\lambda}{\lambda+q_x(a)}p_{n,l}(da|x)\right] \int_{\textbf{A}^I}Q(\Gamma|x,\hat{b}) \overline{p}_{n,k}(d\hat{b}|x)\prod_{m=0}^{k-1}\overline{p}_{n,m}(\textbf{A}^G|x)\\
&=&\sum_{k=0}^\infty  \left[\prod_{l=0}^{k-1} \int_{\textbf{A}^G} \frac{\lambda}{\lambda+q_x(a)} \overline{p}_{n,l}(da|x)\right] \int_{\textbf{A}^I}Q(\Gamma|x,\hat{b}) \overline{p}_{n,k}(d\hat{b}|x)=B_2,
\end{eqnarray*}
where the second to the last equality holds by the definition of $p_{n,l}$: $\overline{p}_{n,l}(da|x)=\overline{p}_{n,l}(\textbf{A}^G|x)p_{n,l}(da|x)$ (see (\ref{2021June30Eqn35})), no matter whether $\prod_{m=0}^{k-1}\overline{p}_{n,m}(\textbf{A}^G|x)$ vanishes or not, and the same remark applies to the calculations for $C_1$ and $C_2$ below, which will not be repeated.

As for $C_1$, we have
\begin{eqnarray*}
&&C_1:=\int_{{\bf \Xi}} p_n(d\xi|x)  \sum_{k=0}^\infty \int_0^{\tau_k} \tilde{q}^\xi(\Gamma|x,t) e^{-\int_0^t q^\xi_x(s)ds}dt \overline{p}_{n,k}(\textbf{A}^I|x)\prod_{m=0}^{k-1}\overline{p}_{n,m}(\textbf{A}^G|x)\\
&=&\sum_{k=0}^\infty \int_{{\bf \Xi}} p_n(d\xi|x)\left( \sum_{l=0}^{k-1}\tilde{q}(\Gamma|x,\alpha_l)\left(\prod_{\nu=0}^{l-1}e^{-\psi_{\nu+1}q_x(\alpha_\nu)}\right) \int_0^{\psi_{l+1}} e^{-t q_x(\alpha_l)}dt \overline{p}_{n,k}(\textbf{A}^I|x)\prod_{m=0}^{k-1}\overline{p}_{n,m}(\textbf{A}^G|x)\right)\\
&=&\sum_{k=0}^\infty \sum_{l=0}^{k-1} \int_{\textbf{A}^G} \frac{\tilde{q}(\Gamma|x,a)}{q_x(a)+\lambda} p_{n,l}(da|x) \left(\prod_{\nu=0}^{l-1} \int_{\textbf{A}^G}\frac{\lambda}{\lambda+q_x(a)}p_{n,\nu}(da|x)\right)\overline{p}_{n,k}(\textbf{A}^I|x)\prod_{m=0}^{k-1}\overline{p}_{n,m}(\textbf{A}^G|x)\\
&=& \sum_{k=0}^\infty \sum_{l=0}^{k-1} \int_{\textbf{A}^G} \frac{\tilde{q}(\Gamma|x,a)}{q_x(a)+\lambda} \overline{p}_{n,l}(da|x) \left(\prod_{\nu=0}^{l-1} \int_{\textbf{A}^G}\frac{\lambda}{\lambda+q_x(a)}\overline{p}_{n,\nu}(da|x)\right)(1-\overline{p}_{n,k}(\textbf{A}^G|x))\prod_{m=l+1}^{k-1}\overline{p}_{n,m}(\textbf{A}^G|x).
\end{eqnarray*}
It is convenient to introduce the following notation:
\begin{eqnarray*}
D_l:=\int_{\textbf{A}^G} \frac{\tilde{q}(\Gamma|x,a)}{q_x(a)+\lambda} \overline{p}_{n,l}(da|x) \left(\prod_{\nu=0}^{l-1} \int_{\textbf{A}^G}\frac{\lambda}{\lambda+q_x(a)}\overline{p}_{n,\nu}(da|x)\right).
\end{eqnarray*}
Then $B_1$ in (\ref{2021June30Eqn12}) can be written as
\begin{eqnarray*}
B_1=\sum_{k=0}^\infty \left(\prod_{i=0}^{k-1}\int_{\textbf{A}^G}\frac{\lambda}{\lambda+q_x(a)}\overline{p}_{n,i}(da|x)\right)\int_{\textbf{A}^G} \frac{\tilde{q}(\Gamma|x,a)}{\lambda+q_x(a)}\overline{p}_{n,k}(da|x)=\sum_{l=0}^\infty D_l,
\end{eqnarray*} which is finite because so is the left hand side of (\ref{2021June30Eqn12}).

With the notation of $D_l,$ now we write
\begin{eqnarray*}
C_1=\sum_{k=0}^\infty \sum_{l=0}^{k-1} D_l(1-\overline{p}_{n,k}(\textbf{A}^G|x))\prod_{m=l+1}^{k-1}\overline{p}_{n,m}(\textbf{A}^G|x).
\end{eqnarray*}
By a similar calculation as for $C_1$, we may write
\begin{eqnarray*}
&&C_2:=\int_{{\bf \Xi}} p_n(d\xi|x) \int_0^\infty \tilde{q}^\xi(\Gamma|x,t)e^{-\int_0^t q^\xi_x(s)ds}dt  \prod_{m=0}^{\infty}\overline{p}_{n,m}(\textbf{A}^G|x)\\
&=&\sum_{k=0}^\infty \int_{\textbf{A}^G} \frac{\tilde{q}(\Gamma|x,a)}{\lambda+q_x(a)}\overline{p}_{n,k}(da|x)\left(\prod_{\nu=0}^{k-1} \int_{\textbf{A}^G}\frac{\lambda}{\lambda+q_x(a)}\overline{p}_{n,\nu}(da|x)\right) \prod_{m\ge k+1}\overline{p}_{n,m}(\textbf{A}^G|x)\\
&=&\sum_{k=0}^{\infty} D_k \prod_{m\ge k+1}\overline{p}_{n,m}(\textbf{A}^G|x).
\end{eqnarray*}
Thus,
\begin{eqnarray*}
&&C_1+C_2=\sum_{k=0}^\infty \sum_{l=0}^{k-1} D_l(1-\overline{p}_{n,k}(\textbf{A}^G|x))\prod_{m=l+1}^{k-1}\overline{p}_{n,m}(\textbf{A}^G|x)+\sum_{k=0}^{\infty} D_k \prod_{m\ge k+1}\overline{p}_{n,m}(\textbf{A}^G|x)\\
&=&\sum_{l=0}^\infty D_l \sum_{k=l+1}^\infty (1-\overline{p}_{n,k}(\textbf{A}^G|x))\prod_{m=l+1}^{k-1}\overline{p}_{n,m}(\textbf{A}^G|x)+\sum_{l=0}^{\infty} D_l \prod_{m\ge l+1}\overline{p}_{n,m}(\textbf{A}^G|x)\\
&=&\sum_{l=0}^\infty D_l \left\{\sum_{k\ge l+1}\left(\prod_{m=l+1}^{k-1}\overline{p}_{n,m}(\textbf{A}^G|x) -\prod_{m=l+1}^{k}\overline{p}_{n,m}(\textbf{A}^G|x)\right)+ \prod_{m\ge l+1}\overline{p}_{n,m}(\textbf{A}^G|x)\right\}\\
&=&\sum_{l=0}^\infty D_l=B_1.
\end{eqnarray*}
(Recall that $\sum_{l=0}^\infty D_l$ converges.)
Combining this with the previous observation, we see that $C_1+C_2+C_3=B_1+B_2,$ and by (\ref{2021June30Eqn12}) and (\ref{2021June30Eqn15}), we see that (\ref{2021June30Eqn10}) holds. Consequently, (\ref{2021June30Eqn16}) follows.

\textit{Step 2.} In view of the definition of $\hat{W}_i(x_0,\sigma^P)$ and $W_i(x_0,\overline{S}^P)$,
it remains to show that
\begin{eqnarray}\label{2021June30Eqn18}
\hat{\rm E}_{x_0}^{\sigma^P}[l_i^{\sigma^P,n}(X_n)]={\rm E}_{x_0}^{\overline{S}^P}\left[I\{X_n\ne x_\infty\} \int_{{\bf \Xi}^{GO}}\int_{(0,\infty]}  \int_0^t c_i^{GO,\xi}(X_n,s)ds {\rm P}_n^{\overline{S}^P,\xi}(\Theta_{n+1}\in dt|X_n)\overline{p}_n(d\xi|X_n)\right].
\end{eqnarray}
for bounded $[0,\infty)$-valued functions $c_i^G,c_i^I$, because the general case can be handled using the monotone convergence theorem.

Note that for each $x\in\textbf{X}$
\begin{eqnarray*}
&&l_i^{\sigma^P,n}(x)= \int_{{\bf \Xi}}\int_{[0,\infty]\times \textbf{A}^I}\int_{0}^\infty   \left\{ \int_0^t c_i^{G,\xi}(x,s)ds+I\{t=\hat{c}<\infty\}c_i^I(x,\hat{b})\right\}\\
&&\times G_n^{\sigma^P,\xi}(dt\times \textbf{X}_\infty|x,\hat{c},\hat{b}) \sigma_n^{P,(0)}(d\hat{c}\times d\hat{b}|x,\xi)p_n(d\xi|x) \\
&=&\int_{{\bf\Xi}}\int_{[0,\infty]\times \textbf{A}^I}\left\{\int_0^{\hat{c}} \int_0^t c_i^{G,\xi}(x,s)ds q_x^\xi(t)e^{-\int_0^t q_x^\xi(s)ds}dt +I\{\hat{c}=\infty\} \int_0^\infty c_i^{G,\xi}(x,s)ds e^{-\int_0^\infty q_x^\xi(s)ds}\right.\\
&&\left.+I\{\hat{c}<\infty\}e^{-\int_0^{\hat{c}} q_x^\xi(s)ds} \left(\int_0^{\hat{c}}c^{G,\xi}_i(x,s)ds+c_i^I(x,\hat{b}) \right)  \right\} \sigma_n^{P,(0)}(d\hat{c}\times d\hat{b}|x,\xi)p_n(d\xi|x)\\
&=&\lim_{m\rightarrow \infty}\int_{{\bf\Xi}}\int_{[0,\infty]\times \textbf{A}^I}\left\{\int_0^{\hat{c}} \int_0^t c^{G,\xi}_i(x,s)e^{-\frac{s}{m}}ds q_x^\xi(t)e^{-\int_0^t q_x^\xi(s)ds}dt\right.\\
&&\left. +I\{\hat{c}=\infty\} \int_0^\infty c_i^{G,\xi}(x,s) e^{-\frac{s}{m}}ds e^{-\int_0^\infty q_x^\xi(s)ds}\right.\\
&&\left.+I\{\hat{c}<\infty\}e^{-\int_0^{\hat{c}} q_x^\xi(s)ds} \left(\int_0^{\hat{c}}c^{G,\xi}_i(x,s)e^{-\frac{s}{m}}ds+c_i^I(x,\hat{b}) \right)  \right\} \sigma_n^{P,(0)}(d\hat{c}\times d\hat{b}|x,\xi)p_n(d\xi|x),
\end{eqnarray*}
where the first equality is by (\ref{2021June30Eqn37}), and the second equality is by (\ref{2021June30Eqn11}).
Applying legitimately integration by parts, we see
\begin{eqnarray*}
\int_0^{\hat{c}} \int_0^t c_i^{G,\xi}(x,s)e^{-\frac{s}{m}}ds q_x^\xi(t)e^{-\int_0^t q_x^\xi(s)ds}dt= \int_0^{\hat{c}} c_i^{G,\xi}(x,t)e^{-\frac{t}{m}}e^{-\int_0^t q_x^\xi(s)ds}dt- e^{-\int_0^{\hat{c}} q_x^\xi(s)ds}\int_0^{\hat{c}}e^{-\frac{s}{m}}c_i^{G,\xi}(x,s)ds,
\end{eqnarray*}
where all the terms are finite, $\hat{c}$ being finite or not, because so are $c_i^I,c_i^G$ assumed. Substituting the previous equality back in the above formula, we see
\begin{eqnarray}\label{2021June30Eqn25}
&&l_i^{\sigma^P,n}(x)=  \lim_{m\rightarrow \infty} \int_{{\bf\Xi}}\int_{[0,\infty]\times \textbf{A}^I}\left\{ \int_0^{\hat{c}} c_i^{G,\xi}(x,t)e^{-\frac{t}{m}}e^{-\int_0^t q_x^\xi(s)ds}dt +I\{\hat{c}<\infty\}e^{-\int_0^{\hat{c}} q_x^\xi(s)ds} c_i^I(x,\hat{b}) \right\}\nonumber \\
&&\times\sigma_n^{P,(0)}(d\hat{c}\times d\hat{b}|x,\xi)p_n(d\xi|x)\nonumber\\
&=&\int_{{\bf\Xi}}\int_{[0,\infty]\times \textbf{A}^I}\left\{ \int_0^{\hat{c}} c_i^{G,\xi}(x,t)e^{-\int_0^t q_x^\xi(s)ds}dt +I\{\hat{c}<\infty\}e^{-\int_0^{\hat{c}} q_x^\xi(s)ds} c_i^I(x,\hat{b}) \right\}\nonumber \\
&&\times \sigma_n^{P,(0)}(d\hat{c}\times d\hat{b}|x,\xi)p_n(d\xi|x),
\end{eqnarray}
where second equality holds by the monotone convergence theorem. Observe that the term inside the parenthesis in the above expression is in the same form as the one in the first equality of (\ref{2021June30Eqn19}), where $\tilde{q}^\xi(\Gamma|x,t)$ and $Q(\Gamma|x,\hat{b})$ are now replaced with $c_i^{G,\xi}(x,t)$ and $c_i^I(x,\hat{b})$, respectively. Therefore, by repeating the calculations  below (\ref{2021June30Eqn19}) in Step 1 with obvious modifications, we see that the following equality holds, which is corresponding to (\ref{2021June30Eqn10}) (or more precisely, the established equality $C_1+C_2+C_3=B_1+B_2$, see more explanations below):
\begin{eqnarray*}
l_i^{\sigma^P,n}(x)=\sum_{k=0}^\infty \left(\prod_{i=0}^{k-1}\int_{\textbf{A}^G}\frac{\lambda}{\lambda+q_x(a)}\overline{p}_{n,i}(da|x)\right)\left(\int_{\textbf{A}^G} \frac{c_i^G(x,a)}{\lambda+q_x(a)}\overline{p}_{n,k}(da|x)+\int_{\textbf{A}^I} c^I_i(x,a)\overline{p}_{n,k}(da|x)\right).
\end{eqnarray*}
Indeed, the term on the right hand side of the above equality corresponds to the term on the right hand side of the second equality in (\ref{2021June30Eqn12}), which coincides with the right hand side of (\ref{2021June30Eqn10}), whereas it was observed earlier that $l_i^{\sigma^P,n}(x)$ corresponds to the left hand side of (\ref{2021June30Eqn10}).

Consequently, the left hand side of (\ref{2021June30Eqn18}) reads
\begin{eqnarray*}
&&\hat{\rm E}_{x_0}^{\sigma^P}[l_i^{\sigma^P,n}(X_n)] =\int_{\textbf{X}}\hat{\rm P}^{\sigma^P}_{x_0}(X_n\in dx)l_i^{\sigma^P,n}(x)= \int_{\textbf{X}}\hat{\rm P}_{x_0}^{\sigma^P}(X_n\in dx)\\
&&\times \left\{\sum_{k=0}^\infty \left(\prod_{i=0}^{k-1}\int_{\textbf{A}^G}\frac{\lambda}{\lambda+q_x(a)}\overline{p}_{n,i}(da|x)\right)\left(\int_{\textbf{A}^G} \frac{c_i^G(x,a)}{\lambda+q_x(a)}\overline{p}_{n,k}(da|x)+\int_{\textbf{A}^I} c^I_i(x,a)\overline{p}_{n,k}(da|x)\right)\right\}.
\end{eqnarray*}

On the other hand,
we may write the right hand side of (\ref{2021June30Eqn18}) as
\begin{eqnarray*}
&& {\rm E}_{x_0}^{\overline{S}^P}\left[I\{X_n\ne x_\infty\} \int_{{\bf \Xi}^{GO}}\int_{(0,\infty]}  \int_0^t c_i^{GO,\xi}(X_n,s)ds {\rm P}_n^{\overline{S}^P,\xi}(\Theta_{n+1}\in dt|X_n)\overline{p}_n(d\xi|X_n)\right]\\
&=& \int_{\textbf{X}}{\rm P}^{\overline{S}^P}_{x_0}(X_n\in dx)   \left\{ \sum_{k=0}^\infty \left(\prod_{i=0}^{k-1}\int_{\textbf{A}}\frac{\overline{\lambda}(a)}{\overline{\lambda}(a)+q^{GO}_x(a)}\overline{p}_{n,i}(da|x)\right)\int_{\textbf{A}} \frac{ c_i^{GO}(x,a)}{\overline{\lambda}(a)+q^{GO}_x(a)}\overline{p}_{n,k}(da|x)\right\}\\
&=&\int_{\textbf{X}}{\rm P}_{x_0}^{\overline{S}^P}(X_n\in dx)\\
&&\left\{\sum_{k=0}^\infty \left(\prod_{i=0}^{k-1}\int_{\textbf{A}^G}\frac{\lambda}{\lambda+q_x(a)}\overline{p}_{n,i}(da|x)\right)\left(\int_{\textbf{A}^G} \frac{c_i^G(x,a)}{\lambda+q_x(a)}\overline{p}_{n,k}(da|x)+\int_{\textbf{A}^I} c^I_i(x,a)\overline{p}_{n,k}(da|x)\right)\right\},
\end{eqnarray*}
where the first equality is by (\ref{2021June30Eqn17}), and the last equality is by the definitions of $\overline{\lambda}(a),$ $\textbf{A}$, $c^{GO}_i$ and $q^{GO}.$
In view of (\ref{2021June30Eqn16}), which was established in the above, we see from the previous equality that (\ref{2021June30Eqn18}) holds, as desired. $\hfill\Box$
\bigskip

\subsection{Proof of Theorem \ref{2021June30Theorem01}}

\par\noindent\textit{Proof of Theorem \ref{2021June30Theorem01}.} In view of the discussions below Proposition \ref{PZ18Theorem01}, we only need show that each policy $\overline{S}$ in the model ${\cal M}^{GO}$ with gradual control only can be replicated by a strategy in the gradual-impulsive control model ${\cal M}.$

According to Theorem 2 of \cite{Piunovskiy:2015} (or Theorem 4.1.1 of \cite{PiunovskiyZhang:2020Book}), for each policy $\overline{S}$ in the model ${\cal M}^{GO}$, there is a replicating Markov policy $\overline{S}^M$ in the same model ${\cal M}^{GO}$ (recall Definition \ref{2021June30Definition01}). Theorem \ref{2021June30Theorem05} and Theorem \ref{2021June30Theorem02} imply that the Markov policy $\overline{S}^M$ in ${\cal M}^{GO}$  is replicated by a Poisson-related strategy $\sigma^P$ in the gradual-impulsive control model ${\cal M}$. To complete the proof of the statement, it remains to show that this replicating Poisson-related strategy $\sigma^P$ in the gradual-impulsive control model ${\cal M}$ can be replicated by an (ordinary) strategy $\sigma$ in the same model ${\cal M}.$ This is justified as follows. Without loss of generality, we assume that $c_i^G$ and $c_i^I$ are nonnegative and bounded in this proof.



Let some Poisson-related strategy $\sigma^P=\{(\sigma^P_n,p_n)\}_{n\ge 0}$ in the model ${\cal M}$ be fixed.

Let
\begin{eqnarray*}
\sigma_n^{(0)}(d\hat{c}\times d\hat{b}|x):=\int_{{\bf\Xi}} \sigma^{P,(0)}_n(d\hat{c}\times d\hat{b}|x,\xi)p_n(d\xi|x).
\end{eqnarray*}
Then, by Proposition 7.27 of \cite{Bertsekas:1978} (or Proposition B.1.33 of \cite{PiunovskiyZhang:2020Book}), there is a stochastic kernel $\hat{p}_n(d\xi|x,\hat{c},\hat{b})$ on ${\cal B}({\bf \Xi})$ given $(x,\hat{c},\hat{b})\in \textbf{X}\times [0,\infty]\times \textbf{A}^I$ satisfying
\begin{eqnarray}\label{2021June30Eqn20}
 \sigma^{P,(0)}_n(d\hat{c}\times d\hat{b}|x,\xi)p_n(d\xi|x)=\hat{p}_n(d\xi|x,\hat{c},\hat{b})\sigma_n^{(0)}(d\hat{c}\times d\hat{b}|x).
\end{eqnarray}

We define a strategy $\sigma={(\sigma_n^{(0)},\hat{F}_n)}_{n\ge 0}$ in the model ${\cal M}$ as follows.
Let \begin{eqnarray*}
\sigma_n^{(0)}(d\hat{c}\times d\hat{b}|\hat{h}_n):=\sigma_n^{(0)}(d\hat{c}\times d\hat{b}|x_n)
\end{eqnarray*}
(Recall the generic notation $\hat{x}_n=(\hat{\theta}_n,x_n)$ for the state in the model ${\cal M}$.)
Let
\begin{eqnarray*}
\hat{F}_n(\hat{h}_n,\hat{c},\hat{b})_t(da):= \frac{\int_{\bf \Xi} e^{-\int_0^t q_{x_n}^\xi(u)du} \sum_{k\ge 0} \delta_{\alpha_k}(da) I\{\tau_k<t\le \tau_{k+1} \}\hat{p}_n(d\xi|x_n,\hat{c},\hat{b})}{\int_{\bf \Xi} e^{-\int_0^t q_{x_n}^\xi(u)du} \hat{p}_n(d\xi|x_n,\hat{c},\hat{b})}=:\hat{F}_n(x_n,\hat{c},\hat{b})(da),
\end{eqnarray*}
where the generic notations $\xi=\{(\psi_n,\alpha_n)\}_{n\ge 0}\in {\bf \Xi}$ and $\tau_k=\sum_{i=0}^k \psi_i$ are in use.

We will show that
\begin{eqnarray}\label{2021June30Eqn21}
\hat{\rm P}_{x_0}^{\sigma}(X_n\in dx)= \hat{\rm P}_{x_0}^{\sigma^P}(X_n\in dx)~\forall~n\ge 0.
\end{eqnarray}
(Recall the generic notation $\hat{X}_n=(\hat{\Theta}_n,X_n)$ in the model ${\cal M}$.) Since the initial states are the same,
with an inductive argument, it is sufficient to show for $\Gamma\in{\cal B}(\textbf{X})$ and $x\in \textbf{X}$,
\begin{eqnarray}\label{2021June30Eqn22}
\hat{\rm P}_{x_0}^{\sigma}(X_{n+1}\in \Gamma|X_n=x)= \hat{\rm P}_{x_0}^{\sigma^P}(X_{n+1}\in \Gamma|X_n=x) ~\forall ~n\ge 0.
\end{eqnarray}

Then,
\begin{eqnarray*}
&&\tilde{q}(\Gamma|x,\hat{F}_n(x,\hat{c},\hat{b})_t)=\int_{\textbf{A}^G} \tilde{q}(\Gamma|x,a) \frac{\int_{{\bf \Xi}}e^{-\int_0^t q_x^\xi(u) du}\sum_{k\ge 0}\delta_{\alpha_k}(da) I\{\tau_k<t\le \tau_{k+1}\} \hat{p}_n(d\xi|x,\hat{c},\hat{b}) }{\int_{\bf \Xi} e^{-\int_0^t q_x^\xi(u)du} \hat{p}_n(d\xi|x,\hat{c},\hat{b})}\\
&=&   \frac{\int_{{\bf \Xi}}\tilde{q}^\xi(\Gamma|x,t)e^{-\int_0^t q_x^\xi(u) du} \hat{p}_n(d\xi|x,\hat{c},\hat{b}) }{\int_{\bf \Xi} e^{-\int_0^t q_x^\xi(u)du} \hat{p}_n(d\xi|x,\hat{c},\hat{b})},
\end{eqnarray*}
recall (\ref{2021June30Eqn32}) for the definition of $\tilde{q}^\xi.$ Applying the above equality to $\Gamma=\textbf{X}$, we see
\begin{eqnarray*}
q_x(\hat{F}_n(x,\hat{c},\hat{b})_t)=\frac{\int_{{\bf \Xi}} q^\xi_x(t)e^{-\int_0^t q_x^\xi(u) du} \hat{p}_n(d\xi|x,\hat{c},\hat{b}) }{\int_{\bf \Xi} e^{-\int_0^t q_x^\xi(u)du} \hat{p}_n(d\xi|x,\hat{c},\hat{b})}=-\frac{d}{dt} \ln \int_{\bf \Xi}   e^{-\int_0^t q_x^\xi(u)du} \hat{p}_n(d\xi|x,\hat{c},\hat{b})
\end{eqnarray*}
for almost all $t$,
and thus
\begin{eqnarray*}
e^{-\int_0^t q_x(\hat{F}_n(x,\hat{c},\hat{b})_s)ds}=\int_{\bf \Xi}   e^{-\int_0^t q_x^\xi(u)du} \hat{p}_n(d\xi|x,\hat{c},\hat{b}).
\end{eqnarray*}

Now,
\begin{eqnarray}\label{2021June30Eqn26}
&&\hat{\rm P}_{x_0}^{\sigma}(X_{n+1}\in \Gamma|X_n=x)=\int_{[0,\infty]\times\textbf{A}^I}\left\{\int_0^{\hat{c}} \tilde{q}(\Gamma|x,\hat{F}_n(x,\hat{c},\hat{b}))_t) e^{-\int_0^ t q_x(\hat{F}_n(x,\hat{c},\hat{b})_s)ds}dt\right.\nonumber\\
&&\left.+I\{\hat{c}<\infty\} e^{-\int_0^{\hat{c}} q_x(\hat{F}_n(x,\hat{c},\hat{b})_s)ds}Q(\Gamma|x,\hat{b})\right\} \sigma_n^{(0)}(d\hat{c}\times d\hat{b}|x)\nonumber\\
&=&\int_{[0,\infty]\times\textbf{A}^I}\left\{\int_0^{\hat{c}} \int_{{\bf \Xi}}\tilde{q}^\xi(\Gamma|x,t)e^{-\int_0^t q_x^\xi(u) du} \hat{p}_n(d\xi|x,\hat{c},\hat{b})dt\right.\nonumber\\
&&\left. +I\{\hat{c}<\infty\} \int_{\bf \Xi}   e^{-\int_0^{\hat{c}} q_x^\xi(u)du} \hat{p}_n(d\xi|x,\hat{c},\hat{b})Q(\Gamma|x,\hat{b}) \right\}\sigma_n^{(0)}(d\hat{c}\times d\hat{b}|x)\nonumber\\
&=&\int_{{\bf \Xi}}  \int_{[0,\infty]\times\textbf{A}^I}\left\{\int_0^{\hat{c}} \tilde{q}^\xi(\Gamma|x,t)e^{-\int_0^t q_x^\xi(u) du} dt\right.\nonumber\\
&&\left. +I\{\hat{c}<\infty\}     e^{-\int_0^{\hat{c}} q_x^\xi(u)du}  Q(\Gamma|x,\hat{b})) \right\}\sigma_n^{P,(0)}(d\hat{c}\times d\hat{b}|x,\xi)p_n(d\xi|x)\nonumber\\
&=&\hat{\rm P}_{x_0}^{\sigma^P}(X_{n+1}\in \Gamma|X_n=x),
\end{eqnarray}
where the second to the last equality is by (\ref{2021June30Eqn20}), and for the last equality, c.f., (\ref{PZ18Eqn12}). Thus, (\ref{2021June30Eqn22}) is verified, and (\ref{2021June30Eqn21}) follows.

Finally, one can show with a similar argument as for (\ref{2021June30Eqn25}) that
\begin{eqnarray*}
&&\hat{\rm E}_{x_0}^{\sigma}\left[l_i(\hat{X}_n,\hat{A}_n,\hat{X}_{n+1})|X_n=x\right]= \int_{[0,\infty]\times \textbf{A}^I}\left\{\int_0^{\hat{c}} c_i^G(x,\hat{F}_n(x,\hat{c},\hat{b})_t)e^{-\int_0^t q_x(\hat{F}_n(x,\hat{c},\hat{b})_s) ds}dt\right.\\
&&\left.+I\{\hat{c}<\infty\}e^{-\int_0^{\hat{c}} q_x(\hat{F}_n(x,\hat{c},\hat{b})_s)ds}c_i^I(x,\hat{b})\right\}\sigma_n^{(0)}(d\hat{c}\times d\hat{b}|x),
\end{eqnarray*}
where $l_i$ was defined by (\ref{PZ18Eqn11}).
Having inspected that the term in the parenthesis of the last equality is in the same form as the term on the right hand side of the first equality in  (\ref{2021June30Eqn26}), we see now
\begin{eqnarray*}
&&\hat{\rm E}_{x_0}^{\sigma}\left[l_i(\hat{X}_n,\hat{A}_n,\hat{X}_{n+1})|X_n=x\right]=\int_{{\bf \Xi}}  \int_{[0,\infty]\times\textbf{A}^I}\left\{\int_0^{\hat{c}} c_i^{G,\xi}(x,t)e^{-\int_0^t q_x^\xi(u) du} dt\right.\nonumber\\
&&\left. +I\{\hat{c}<\infty\}     e^{-\int_0^{\hat{c}} q_x^\xi(u)du}  c^I_i(x,\hat{b})) \right\}\sigma_n^{P,(0)}(d\hat{c}\times d\hat{b}|x,\xi)p_n(d\xi|x)= l_i^{\sigma^P,n}(x),
\end{eqnarray*}
where the first equality corresponds to the second to the last equality in (\ref{2021June30Eqn26}), and the last equality holds by (\ref{2021June30Eqn25}). The previous equality and (\ref{2021June30Eqn21}) imply that
\begin{eqnarray*}
\hat{\rm E}_{x_0}^{\sigma}\left[l_i(\hat{X}_n,\hat{A}_n,\hat{X}_{n+1})\right]=\hat{\rm E}_{x_0}^{\sigma^P}\left[l_i^{\sigma^P,n}(X_n)\right]
\end{eqnarray*}
for all $n\ge 0$. The statement is thus proved. $\hfill\Box$
\bigskip

\section{Conclusion}\label{2021June30Sec06}
In conclusion, we investigated a constrained optimal control problem for a gradual-impulsive CTMDP with the performance criteria being the total undiscounted costs. We fully justified a reduction method, and thus closed an open issue left in \cite{PiunovskiyZhang:2021SICON}. The reduction method induces an equivalent but simpler standard CTMDP model. The effectiveness of this method was demonstrated when we used it to establish, under rather natural conditions, the linear programming approach to solving the concerned constrained optimal control problem.

\appendix

\section{Appendix: Proof of Proposition \ref{2021June30Proposition01}}\label{2021June30Subsection01}

\par\noindent\textit{Proof of Proposition \ref{2021June30Proposition01}.} Since the initial states $X_0$ are the same in both models ${\cal M}, {\cal M}^{GO}$, it is sufficient to show the following: if the current state is $X_n=x$, then the distribution of the state $X_{n+1}$ after the next sojourn time, as well as the expected accumulated cost over the next sojourn time in ${\cal M}$ under $\sigma^S=(\sigma^{S,(0)},\hat{F}^S)$ and in ${\cal M}^{GO}$ under $\overline{F}^S$ coincide. We will verify this for $x\in O$ and $x\in \textbf{X}\setminus O$: the case when $x=x_\infty$ is trivial with the next state being $x_\infty$ and the accumulated cost being $0$.

Suppose $x\in O$. Note that
\begin{eqnarray*}
O=\left\{x\in \textbf{X}:~\int_{\textbf{A}}q^{GO}_{x}(a)\overline{F}^S(x)(da)>0\right\}=\left\{x\in \textbf{X}:~ \int_{\textbf{A}^G} q_x(a) \overline{F}^S(x)(da) +  \overline{F}^S(x)(\textbf{A}^I)>0\right\}
\end{eqnarray*}
by the definition of $q_x^{GO}(a)$ and $\textbf{A},$ which will be used below in this proof without special reference.  Then in the model ${\cal M}^{GO}$ under the stationary policy $\overline{F}^S$,
\begin{eqnarray*}
&&{\rm P}^{\overline{F}^S}_{x_0}(X_{n+1}\in dy|X_n=x)= \frac{\int_{\textbf{A}}\tilde{q}^{GO}(dy|x,a)\overline{F}^S(x)(da)}{\int_{\textbf{A}}{q}^{GO}_x(a)\overline{F}^S(x)(da)}\\
&=&\frac{\int_{\textbf{A}^G}\tilde{q}(dy|x,a)\overline{F}^S(x)(da)+\int_{\textbf{A}^I}Q(dy|x,b)\overline{F}^S(x)(db)}{\int_{\textbf{A}^G}{q}_x(a)\overline{F}^S(x)(da)+\overline{F}^S(\textbf{A}^I)}~\mbox{on ${\cal B}(\textbf{X})$}.
\end{eqnarray*}
(The denominator is not vanishing because $x\in O.$) On the other hand, in the model ${\cal M}$ under the stationary strategy $\sigma^S=(\sigma^{S,(0)},\hat{F}^S)$, we consider three cases
\begin{itemize}
\item[(a)] $\int_{\textbf{A}^G}q_x(a)\overline{F}^S(x)(da)>0$ and $\overline{F}^S(x)(\textbf{A}^I)>0$,
\item[(b)]  $\int_{\textbf{A}^G}q_x(a)\overline{F}^S(x)(da)>0$ and $\overline{F}^S(x)(\textbf{A}^I)=0$,
\item[(c)]  $\int_{\textbf{A}^G}q_x(a)\overline{F}^S(x)(da)=0$.
\end{itemize}

Case (a): if $\int_{\textbf{A}^G}q_x(a)\overline{F}^S(x)(da)>0$ and $\overline{F}^S(x)(\textbf{A}^I)>0$, then $\overline{F}^S(x)(\textbf{A}^G)>0$, and
\begin{eqnarray*}
&&\hat{{\rm P}}_{x_0}^{\sigma^S}(X_{n+1}\in dy|X_n=x)\\
&=&\sigma^{S,(0)}(\{\infty\}\times \textbf{A}^I|x) \frac{\int_{\textbf{A}^G} \tilde{q}(dy|x,a)\hat{F}^S(x)(da)}{\int_{\textbf{A}^G} q_x(a)\hat{F}^S(x)(da)}+\int_{\textbf{A}^I}Q(dy|x,\hat{b})\sigma^{S,(0)}(\{0\}\times d\hat{b}|x)\\
&=& \frac{\int_{\textbf{A}^G} q_x(a) \overline{F}^S(x)(da) }{\int_{\textbf{A}^G} q_x(a) \overline{F}^S(x)(da) +  \overline{F}^S(x)(\textbf{A}^I)} \frac{\int_{\textbf{A}^G} \tilde{q}(dy|x,a)\overline{F}^S(x)(da)}{\int_{\textbf{A}^G} q_x(a)\overline{F}^S(x)(da)}\\
&&+\int_{\textbf{A}^I}Q(dy|x,\hat{b})\frac{\overline{F}^S(x)(d\hat{b})}{\int_{\textbf{A}^G} q_x(a) \overline{F}^S(x)(da) +  \overline{F}^S(x)(\textbf{A}^I)}\\
&=&\frac{\int_{\textbf{A}^G}\tilde{q}(dy|x,a)\overline{F}^S(x)(da)+\int_{\textbf{A}^I}Q(dy|x,b)\overline{F}^S(x)(db)}{\int_{\textbf{A}^G} q_x(a) \overline{F}^S(x)(da) +  \overline{F}^S(x)(\textbf{A}^I)}={\rm P}^{\overline{F}^S}_{x_0}(X_{n+1}\in dy|X_n=x) ~\mbox{on ${\cal B}(\textbf{X})$}.
\end{eqnarray*}

Case (b): if $\int_{\textbf{A}^G}q_x(a)\overline{F}^S(x)(da)>0$ and $\overline{F}^S(x)(\textbf{A}^I)=0$, then $\overline{F}^S(x)(\textbf{A}^G)>0$, and $\sigma^{S,(0)}(\{\infty\}\times \textbf{A}^I|x)=1,$  so that
\begin{eqnarray*}
&&\hat{{\rm P}}_{x_0}^{\sigma^S}(X_{n+1}\in dy|X_n=x)=\sigma^{S,(0)}(\{\infty\}\times \textbf{A}^I|x) \frac{\int_{\textbf{A}^G} \tilde{q}(dy|x,a)\hat{F}^S(x)(da)}{\int_{\textbf{A}^G} q_x(a)\hat{F}^S(x)(da)}\\
&=& \frac{\int_{\textbf{A}^G} q_x(a) \overline{F}^S(x)(da) }{\int_{\textbf{A}^G} q_x(a) \overline{F}^S(x)(da)  } \frac{\int_{\textbf{A}^G} \tilde{q}(dy|x,a)\overline{F}^S(x)(da)}{\int_{\textbf{A}^G} q_x(a)\overline{F}^S(x)(da)}\\
&=&\frac{\int_{\textbf{A}^G}\tilde{q}(dy|x,a)\overline{F}^S(x)(da)}{\int_{\textbf{A}^G} q_x(a) \overline{F}^S(x)(da) }=\frac{\int_{\textbf{A}^G}\tilde{q}(dy|x,a)\overline{F}^S(x)(da)+\int_{\textbf{A}^I}Q(dy|x,b)\overline{F}^S(x)(db)}{\int_{\textbf{A}^G}{q}_x(a)\overline{F}^S(x)(da)+\overline{F}^S(x)(\textbf{A}^I)}\\ &=&  \frac{\int_{\textbf{A}} \tilde{q}^{GO}(dy|x,a)\overline{F}^S(x)(da) }{\int_{\textbf{A}} q^{GO}_x(a)\overline{F}^S(x)(da) }  ={\rm P}^{\overline{F}^S}_{x_0}(X_{n+1}\in dy|X_n=x)~\mbox{on ${\cal B}(\textbf{X})$}.
\end{eqnarray*}
where the second to the last equality holds because $\overline{F}^S(x)(\textbf{A}^I)=0$.

Case (c): if  $\int_{\textbf{A}^G}q_x(a)\overline{F}^S(x)(da)=0$, then $\overline{F}^S(x)(\textbf{A}^I)>0$ (since $x\in O$),  $\sigma^{S,(0)}(\{0\}\times \textbf{A}^I|x)=1,$ and  $\sigma^{S,(0)}(\{0\}\times d\hat{b}|x)=\frac{\overline{F}^S(x)(d\hat{b})}{\overline{F}^S(x)(\textbf{A}^I)}$ in which case,
 \begin{eqnarray*}
 &&\hat{{\rm P}}_{x_0}^{\sigma^S}(X_{n+1}\in dy|X_n=x)= \int_{\textbf{A}^I}Q(dy|x,\hat{b})\sigma^{S,(0)}(\{0\}\times d\hat{b}|x)=\int_{\textbf{A}^I}Q(dy|x,\hat{b}) \frac{\overline{F}^S(x)(d\hat{b})}{\overline{F}^S(x)(\textbf{A}^I)} \\
 &=&\frac{\int_{\textbf{A}^G}\tilde{q}(dy|x,a)\overline{F}^S(x)(da)+\int_{\textbf{A}^I}Q(dy|x,b)\overline{F}^S(x)(db)}{\int_{\textbf{A}^G}{q}_x(a)\overline{F}^S(x)(da)+\overline{F}^S(x)(\textbf{A}^I)}={\rm P}^{\overline{F}^S}_{x_0}(X_{n+1}\in dy|X_n=x)~\mbox{on ${\cal B}(\textbf{X}).$}
 \end{eqnarray*}

Now suppose $x\in \textbf{X}\setminus O,$ that is,
\begin{eqnarray*}
\int_{\textbf{A}}q^{GO}_x(a)\overline{F}^S(x)(da)= \int_{\textbf{A}^G} q_x(a) \overline{F}^S(x)(da) +  \overline{F}^S(x)(\textbf{A}^I)=0,
\end{eqnarray*}
and in particular, $\overline{F}^S(x)(\textbf{A}^I)=0$, $\overline{F}^S(x)(\textbf{A}^G)=1$, and $\sigma^{S,(0)}(\{\infty\}\times \textbf{A}^I|x)=1.$ Then
\begin{eqnarray*}
\hat{{\rm P}}_{x_0}^{\sigma^S}(X_{n+1}=x_\infty|X_n=x)={\rm P}^{\overline{F}^S}_{x_0}(X_{n+1}=x_\infty|X_n=x)=1.
\end{eqnarray*}

Thus, we have verified
\begin{eqnarray*}
\hat{{\rm P}}_{x_0}^{\sigma^S}(X_{n+1}\in dy|X_n=x)={\rm P}^{\overline{F}^S}_{x_0}(X_{n+1}\in dy|X_n=x).
\end{eqnarray*}

The similar argument can be used to show that when $c_i^G$ and $c_i^I$ are $[0,\infty]$-valued, given $X_n=x,$ the expected accumulated costs over the next sojourn time in both models ${\cal M},{\cal M}^{GO}$, under $\sigma^S$ and $\overline{F}^S$ respectively, are both given by
\begin{eqnarray*}
\frac{\int_{\textbf{A}}c_i^{GO}(x,a)\overline{F}^S(x)(da)}{\int_{\textbf{A}}q_x^{GO}(a)\overline{F}^S(x)(da)}=
\frac{\int_{\textbf{A}^G}c_i^G(x,a)\overline{F}^S(x)(da)+ \int_{\textbf{A}^I}c_i^I(x,b)\overline{F}^S(x)(db)}{\int_{\textbf{A}^G}q_x(a)\overline{F}^S(x)(da)+ \overline{F}^S(x)(\textbf{A}^I)},
\end{eqnarray*}
where it is accepted that $\frac{0}{0}:=0.$ The case of general-signed cost rate and function follows from this by considering the positive and negative parts. The proof of the statement is now complete. $\hfill\Box$
\bigskip




\subsection*{Acknowledgement} This paper was discussed at the Liverpool workshop: modern
trends in controlled stochastic processes (July, 2021), supported by the EPSRC (EP/T018216/1).

\end{document}